\documentclass[12pt]{article}
\usepackage[cp866]{inputenc}
\usepackage[psamsfonts]{amssymb}
\usepackage{amsmath,amsfonts,amssymb,amscd}
\usepackage{graphicx}
\begin{document}

\title{Representations of virtual braids by automorphisms and virtual knot groups}

\author{{V.~G.~Bardakov, Yu. A. Mikhalchishina, M.~V.~Neshchadim }}%

\maketitle {\small}

\vspace{0.8cm}


\begin{abstract}
In the present paper the representation of the virtual braid group $VB_n$ into the automorphism group of free product of the free group and free abelian group is constructed.
This representation generalizes the previously constructed ones. The fact that these already known representations are not faithful for $n \geq 4$ is verified. Using representations of $VB_n$,
the virtual link group is defined. Also representations of welded braid group $WB_n$ are constructed and the welded link group is defined.
\end{abstract}

\begin{center}
{\bf Introduction }
\end{center}

The Artin representation \cite[Corollary 1.8.3]{Bir} of the braid group by automorphisms of the free group allows to solve the word problem in the braid group, to construct the linear representations of the braid group (the Burau representation) as well as the linear representation of pure braid group (the Gassner representation) etc. The virtual braid group  (see \cite{V}, \cite{Ka}) is one of generalizations of the classical braid group.
A natural question about the representation of the virtual braid group by automorphisms of some group arises.

For each virtual braid the virtual link could be constructed. Contrary to the classical link case there is no natural definition of a virtual link group.
The following question arises: to define the virtual link group.

There exist several approaches to determine the group of virtual knots and links.
The first definition was given by L. Kauffman \cite{Ka} in 1996. Next there was the definition by
V. Manturov \cite{Man} (see also \cite{BB}, where the same groups were defined using the virtual braid group representation by automorphisms of a free group).
Further there was the definition by D. Silver and S. Williams \cite{SW}-\cite{SW1} and lately another version by
H. Boden,  E. Dies, A. Gaudreau,  A. Gerlings,  E. Harper,  A. Nicas ~\cite{BDG}.

The fact that all representations known before are not faithful for $n \geq 4$ will be shown in the present work. A new representation which generalizes all previous ones will be suggested. Moreover the fact that the new representation is equivalent to the extension of the Artin representation of the braid group will be shown. Using the new representation the virtual link group is defined and the fact that it is the virtual link invariant is proved.
There will be shown the fact that for the classical link the constructed group is isomorphic to the free product of the fundamental group of the link complement in the tree dimensional sphere and some free abelian group.

The principally new way to the definition of the virtual link group was suggested in the paper of J. Carter, D. Silver and S. Williams \cite{CSW}
where the virtual link group was defined as the fundamental group of the complement of some classical link in a thickened surface.
In the present work there will be shown that the group of virtual trefoil  defined in  \cite{CSW} is not isomorphic to the group of virtual trefoil  defined using the representation of virtual braid group .

Authors are grateful to Oleg Chterental who pointed out an element in the kernel  of the Artin representation of the virtual braid group (see \cite{Ch}).

\vspace{0.8cm}


\begin{center}
{\bf \S~1. The braid group, the virtual braid group and the welded braid group}
\end{center}

\vspace{0.5cm}

In the paper the following notations are used. If  $a$ and $b$ are elements of some group $G$, then $a^b = b^{-1} a b$ is conjugacy of  $a$ by  $b$;
$[a, b]=a^{-1}b^{-1} a b$ is the commutator of elements~$a$ and~$b$.

{\bf The braid group.} Recall that the braid group $B_{n}$  on $n$ strands, where $n\geq 2$, is defined by generators $\sigma_{1},\sigma_{2}, \ldots,\sigma_{n-1}$ and defining relations
 \begin{align}
 \sigma_{i}\sigma_{i+1}\sigma_{i} = ~&\sigma_{i+1}\sigma_{i}\sigma_{i+1} \hspace{2cm} \mbox{for} ~ i = 1, 2, \ldots,n-2,\label{e:CIF1}\\
  \sigma_{i}\sigma_{j} = ~&\sigma_{j}\sigma_{i} \hspace{3cm} ~\mbox{for}~ |i-j|\geq 2.\label{e:CIF2}
\end{align}
There exists a homomorphism of the group $B_{n}$ on the permutation group
$S_{n}$ that maps the generator $\sigma_{i}$ into
 the transposition $(i,{i+1}),$ $ i = 1, 2, \ldots,n-1.$ The  kernel of this homomorphism is {\it the pure braid group} and is denoted by the symbol $P_{n}$. The group $P_{n}$ is generated by
elements $a_{ij}, 1\leq i< j\leq n$, which can be expressed through the  generators of $B_{n}$ in the following way
\[
 a_{i,i+1} = \sigma_{i}^{2},
\]
\[a_{ij} =
 \sigma_{j-1}\sigma_{j-2} \ldots\sigma_{i+1}\sigma_{i}^{2}\sigma_{i+1}^{-1} \ldots
 \sigma_{j}^{-1}\sigma_{j-1}^{-1}.
\]
The group $P_n$ composes in the semi direct product of free groups
$$
P_n=U_n\leftthreetimes (U_{n-1}\leftthreetimes (\ldots \leftthreetimes
(U_3\leftthreetimes U_2))\ldots),
$$
where $U_i = \langle a_{1i}, a_{2i}, \ldots, a_{i-1,i} \rangle \simeq F_{i-1}$ is the free group of the rank $n-1$, $i=2,3,\ldots,n$ (see \cite[Section 1.4]{Bir},  \cite{Mar}).

{\bf The virtual braid group and the welded braid group.}
{\it The virtual braid group} $VB_{n}$ is presented in the paper \cite{Ka}, where a
 system of its generators and defining relations is written. In the paper \cite{V}
the more compact system of defining relations is constructed (see relations below).
The virtual braid group $VB_{n}$ is generated by the classical braid group
$B_n = \langle \sigma_1, \ldots,\sigma_{n-1} \rangle$
and the permutation group $S_n = \langle\rho_1, \ldots,\rho_{n-1}\rangle$.
Generators $\sigma_{i}, i = 1, \ldots, n-1$
satisfy the relations (\ref{e:CIF1}), (\ref{e:CIF2}) and generators $\rho_{i}, i = 1, \ldots, n-1,$  satisfy the following relations (relations in the permutation group $S_n$):
\begin{align}
  \rho_{i}^{2} &= 1 ~~ \hspace{3,5cm} \mbox{for}~~ i = 1, 2, \ldots,n-1,\label{e:CIF8}\\
 \rho_{i}\rho_{j} &= \rho_{j}\rho_{i} ~~      \hspace{3cm} \mbox{for}~~ |i-j|\geq 2,\label{e:CIF9}\\
 \rho_{i}\rho_{i+1}\rho_{i} &= \rho_{i+1}\rho_{i}\rho_{i+1}~~ \hspace{2cm} \mbox{for}~~ i = 1, 2 \ldots,n-2.\label{e:CIF10}\\
 \intertext{Other defining relations of the group  $VB_{n}$ are mixed and they are as follows}
 \sigma_{i}\rho_{j}& = \rho_{j}\sigma_{i} ~~      \hspace{3cm} \mbox{for}~~ |i-j|\geq 2,\label{e:CIF11}\\
 \rho_{i}\rho_{i+1}\sigma_{i} &= \sigma_{i+1}\rho_{i}\rho_{i+1}~~ \hspace{2cm} \mbox{for}~~ i = 1, 2, \ldots,n-2.\label{e:CIF12}\\
 \intertext{Note that the last relation is equivalent to:}
 \rho_{i+1}\rho_{i}\sigma_{i+1} &= \sigma_{i}\rho_{i+1}\rho_{i}~~ \hspace{2,5cm} \mbox{for}~~ i = 1, 2, \ldots,n-2.\label{e:CIF13}
 \end{align}
 As it was pointed out in the paper \cite{Ka}, in the group $VB_{n}$ the following relations do not hold
\begin{align*}
 \mathcal{F}_1: ~~~~\rho_{i}\sigma_{i+1}\sigma_{i} &
      =\sigma_{i+1}\sigma_{i}\rho_{i+1}, &\mbox{for}~~ i = 1, 2, \ldots,n-2,\\
 \mathcal{F}_2: ~\rho_{i+1}\sigma_{i}\sigma_{i+1} &
      =\sigma_{i}\sigma_{i+1}\rho_{i}, &\mbox{for}~~ i = 1, 2, \ldots,n-2,
\end{align*}
which look like relations (\ref{e:CIF13}).
These relations are referred to as forbidden relations.

In the work \cite{F} the welded braid group denoted by $WB_{n}$ is presented. This group is generated by elements $\sigma_{i}$ and $\alpha_{i},~ i = 1,
2, \ldots, n-1$. The group generated by elements $\sigma_{i}$ is the classical braid group $B_{n}$, the group generated by elements $\alpha_{i}$ is the permutation group $S_{n}$ also the mixed relations hold
\begin{align}
 \alpha_{i}\sigma_{j} &= \sigma_{j}\alpha_{i}       \hspace{3,5cm} ~\mbox{for}~~ |i-j|\geq 2, \label{e:CIF14}\\
 \alpha_{i}\alpha_{i+1}\sigma_{i} &=\sigma_{i+1}\alpha_{i}\alpha_{i+1}  \hspace{2,5cm} \mbox{for}~~ i = 1, 2,\ldots,n-2,\label{e:CIF15}\\
  \sigma_{i}\sigma_{i+1}\alpha_{i} &= \alpha_{i+1}\sigma_{i}\sigma_{i+1} \hspace{2,5cm} \mbox{for}~~ i = 1, 2, \ldots, n-2.\label{e:CIF16}
 \end{align}
Comparing the relations of groups $VB_{n}$ and $WB_{n}$ it is clear that
$WB_{n}$ is obtained from $VB_{n}$ by adding the additional relation (\ref{e:CIF16}) which coincide with $\mathcal{F}_2$.
Therefore, there exists the homomorphism
$$\varphi_{VW}
: VB_{n}\rightarrow WB_{n},
$$
putting $\sigma_{i}$ into $\sigma_{i}$ and $\rho_{i}$ into
$\alpha_{i}$ for all $i = 1, 2,\ldots,n-1$. Hence $WB_{n}$ is the homomorphic image of the group $VB_{n}$.\\

Note that $WB_n$ is isomorphic to the group of conjugating automorphisms $C_n$, which is a subgroup of $\mathrm{Aut}(F_n)$ and is studied in \cite{B, Bar, Sav, Sav1}.

\vspace{0.8cm}

\begin{center}
{\bf \S~2. Representations of virtual braids by automorphisms}
\end{center}

\vspace{0.5cm}

The present section we begin by listing some known representations of the virtual braid group $VB_n$ by automorphisms of the free group. Next we define the representation that generalizes all previous ones and in the conclusion we show that all previously written representations are not faithful for $n \geq 4$.

{\bf The extension of the Artin representation.}
The representation $VB_n \longrightarrow \mathrm{Aut}(F_{n+1})$ (see \cite{Bar-1}, \cite{Man})
into the automorphism group of the free group $F_{n+1} = \langle x_1, x_2, \ldots, x_n, y \rangle$ of the rank $n+1$:
$$
\varphi_A(\sigma_i) :
\left\{
\begin{array}{ll}
  x_i \longmapsto x_i x_{i+1} x_i^{-1}, &\\
  x_{i+1} \longmapsto x_i, & \\
\end{array}
\right.~~~~~
\varphi_A(\rho_i) :
\left\{
\begin{array}{ll}
  x_i \longmapsto  x_{i+1}^{y^{-1}}, &\\
  x_{i+1} \longmapsto x_i^y, & \\
\end{array}
\right.
$$
(here and hereafter we write only non-trivial actions on the generators, assume that all other generators are fixed) is an extension of the Artin representation
$\varphi_A : B_n \longrightarrow \mathrm{Aut}(F_{n})$ (see \cite{Bir}).
Hence we will hereafter denote it with the same symbol $\varphi_A$.

\bigskip

{\bf The Silver-Williams representation.}
For the pair of natural numbers $n$ and $k$ we define the group
$F_{n,k+1} = F_n * \mathbb{Z}^{k+1}$ as the free product of the free group of the rank $n$ and the free abelian group of the rank $k+1$. We assume that $F_n$ is freely generated by elements
$x_1, x_2, \ldots, x_n$ and $\mathbb{Z}^{k+1}$ is freely generated by elements $v, u_1, u_2, \ldots, u_k$.
Note that $F_{n,1} = F_{n+1}$.

Using the definition of the generalized Alexander group for virtual links given by D. Silver and S. Williams \cite{SW}, we construct the representation
$\varphi_{SW} : VB_n \longrightarrow \mathrm{Aut}(F_{n,n+1})$ of the virtual braid group $VB_n$ into the automorphism group of the group $F_{n,n+1}$ defined by the action on the generators
$$
\varphi_{SW}(\sigma_i) :
\left\{
\begin{array}{ll}
  x_i \longmapsto  x_i x_{i+1}^{u_i} x_i^{-v u_{i+1}}, &\\
  x_{i+1} \longmapsto x_i^v, & \\
\end{array}
\right.~~~
\varphi_{SW}(\sigma_i) :
\left\{
\begin{array}{ll}
  u_i \longmapsto  u_{i+1}, &\\
  u_{i+1} \longmapsto u_i, & \\
\end{array}
\right.~~~
$$
$$
\varphi_{SW}(\rho_i) :
\left\{
\begin{array}{ll}
  x_i \longmapsto   x_{i+1}, &\\
  x_{i+1} \longmapsto x_i, & \\
\end{array}
\right.~~~
\varphi_{SW}(\rho_i) :
\left\{
\begin{array}{ll}
  u_i \longmapsto  u_{i+1}, &\\
  u_{i+1} \longmapsto u_i, & \\
\end{array}
\right.~~~
$$

It is easy to verify that the following automorphism corresponds to the element $\sigma_i^{-1}$
$$
\varphi_{SW}(\sigma_i^{-1}) :
\left\{
\begin{array}{ll}
  x_{i} \longmapsto x_{i+1}^{v^{-1}}, & \\
  x_{i+1} \longmapsto  \left( x_{i+1}^{-v^{-1}} x_{i} x_{i+1}^{u_{i}} \right)^{u_{i+1}^{-1}}, &\\
\end{array}
\right.~~~
\varphi_{SW}(\sigma_i^{-1}) :
\left\{
\begin{array}{ll}
  u_i \longmapsto  u_{i+1}, &\\
  u_{i+1} \longmapsto u_i. & \\
\end{array}
\right.~~~
$$

We shall prove that this representation generalizes the representation
$\varphi_A : VB_n \longrightarrow \mathrm{Aut}(F_{n+1})$.
More precisely the following proposition is true (The idea of this proposition was suggested by D. Silver).

\medskip

{\bf Proposition 1.} {\it If we put
$$
u_1 = u_2 = \ldots = u_n = y, v= y^{-1},
$$
and instead of generators $x_1, x_2, \ldots, x_n$ we put new generators
$$
z_1 = x_1, z_2 = x_2^y, z_3 = x_3^{y^2}, \ldots, z_n = x_n^{y^{n-1}},
$$
then the obtained group is isomorphic to  $F_{n+1}$ and the representation
$\varphi_{SW}$ is coincide with the representation $\varphi_A$ on this group.}

\medskip

{\bf Proof.}
The Silver-Williams representation acts on generators $z_1, z_2, \ldots, z_n$ by the rule
$$
\varphi_{SW}(\sigma_i) :
\left\{
\begin{array}{ll}
  z_i \longmapsto  z_i z_{i+1} z_i^{-1}, &\\
  z_{i+1} \longmapsto z_i, & \\
\end{array}
\right.
\varphi_{SW}(\rho_i) :
\left\{
\begin{array}{ll}
  z_i \longmapsto  z_{i+1}^{y^{-1}}, &\\
  z_{i+1} \longmapsto z_i^y, & \\
\end{array}
\right.
$$
but it is the representation $\varphi_A$.

\bigskip

{\bf The Boden - Dies - Gaudreau - Gerlings - Harper - Nicas representation} \cite{BDG}.
For brevity sake we  denote this representation by the symbol $\varphi_{BD}$.
Let  $F_{n,2} = F_n * \mathbb{Z}^{2}$, where
 $F_n$ is freely generated by elements
$x_1, x_2, \ldots, x_n$ and $\mathbb{Z}^{2}$ is freely generated by elements $v, u$.
The representation
$\varphi_{BD} : VB_n \longrightarrow \mathrm{Aut}(F_{n,2})$
of the virtual braid group $VB_n$ into the automorphism group of the group $F_{n,2}$ is defined by the action on the generators:
$$
\varphi_{BD}(\sigma_i) :
\left\{
\begin{array}{ll}
  x_i \longmapsto  x_i x_{i+1} x_i^{-u}, &\\
  x_{i+1} \longmapsto x_i^u, & \\
\end{array}
\right.~~~
\varphi_{BD}(\rho_i) :
\left\{
\begin{array}{ll}
  x_i \longmapsto   x_{i+1}^{v^{-1}}, &\\
  x_{i+1} \longmapsto x_i^v. & \\
\end{array}
\right.
$$

{\bf A new representation of the virtual braid group.}
Now we define a new representation of $VB_n$ which generalizes the representations listed above.
We consider the free product $F_{n,2n+1} = F_n * \mathbb{Z}^{2n+1}$, where $F_n$ is a free group of the rank $n$ generated by elements
$x_1, x_2, \ldots, x_n$ and $\mathbb{Z}^{2n+1}$ is a free abelian group of the rank $2n+1$ freely generated by elements $u_1, u_2, \ldots, u_n$, $v_0, v_1, v_2, \ldots, v_n$. The following statement is true

\medskip

{\bf Theorem 1.} {\it The following mapping
$\varphi_{M} : VB_n \longrightarrow \mathrm{Aut}(F_{n,2n+1})$ defined by the action on the generators:
$$
\varphi_{M}(\sigma_i) :
\left\{
\begin{array}{l}
  x_i \longmapsto  x_i x_{i+1}^{u_i} x_i^{- v_0 u_{i+1}}, \\
  x_{i+1} \longmapsto x_i^{v_0},  \\
\end{array}
\right.~~~
\varphi_{M}(\sigma_i) :
\left\{
\begin{array}{l}
  u_i \longmapsto  u_{i+1}, \\
  u_{i+1} \longmapsto u_i,  \\
\end{array}
\right.
$$
$$
\varphi_{M}(\sigma_i) :
\left\{
\begin{array}{l}
  v_i \longmapsto  v_{i+1}, \\
  v_{i+1} \longmapsto v_i,  \\
\end{array}
\right.
$$
$$
\varphi_{M}(\rho_i) :
\left\{
\begin{array}{l}
  x_i \longmapsto  x_{i+1}^{v_i^{-1}}, \\
  x_{i+1} \longmapsto x_i^{v_{i+1}},  \\
\end{array}
\right.~~~
\varphi_{M}(\rho_i) :
\left\{
\begin{array}{l}
  u_i \longmapsto  u_{i+1}, \\
  u_{i+1} \longmapsto u_i,  \\
\end{array}
\right.
$$
$$
\varphi_{M}(\rho_i) :
\left\{
\begin{array}{l}
  v_i \longmapsto  v_{i+1}, \\
  v_{i+1} \longmapsto v_i,  \\
\end{array}
\right.
$$
is provided a representation of  $VB_n$ into } $\mathrm{Aut}(F_{n,2n+1})$.

{\bf Proof.}  It is need to verify that the relations of the group $VB_n$
go to the relations of $VB_n$. It is easy to check that it is so for the commutativity relations since the representation $\varphi_M$ is the local representation.

Let us verify the  correctness of the relation
$$
\varphi_M (\sigma_i) \varphi_M (\sigma_{i+1}) \varphi_M (\sigma_i) =
 \varphi_M (\sigma_{i+1}) \varphi_M (\sigma_i) \varphi_M(\sigma_{i+1}).
$$
For that purpose we find the  action of
$\varphi_M (\sigma_i) \varphi_M (\sigma_{i+1}) \varphi_M(\sigma_i)$
on the generators $x_1$, $\ldots$, $x_n$. We have
$$
\varphi_M (\sigma_i) \varphi_M (\sigma_{i+1}) \varphi_M (\sigma_i) :
\left\{
\begin{array}{ll}
  x_{i} \longmapsto x_i  x_{i+1}^{u_i} x_{i+2}^{u_i u_{i+1}} x_{i+1}^{-u_{i} u_{i+2} v_0} x_i^{-u_{i+2} v_0}, & \\
  x_{i+1} \longmapsto  {(x_i x_{i+1}^{u_i} x_{i}^{-u_{i+1} v_0})}^{v_0}, &\\
  x_{i+2} \longmapsto  x_i^{v_0^2}. &\\
\end{array}
\right.
$$
Analogously we find  the  action of
$\varphi_M (\sigma_{i+1}) \varphi_M (\sigma_{i}) \varphi_M (\sigma_{i+1})$. We obtain
$$
\varphi_M (\sigma_{i+1}) \varphi_M (\sigma_{i}) \varphi_M (\sigma_{i+1}) :
\left\{
\begin{array}{ll}
  x_{i} \longmapsto x_i  \left( x_{i+1} x_{i+2}^{u_{i+1}} x_{i+1}^{-u_{i+2} v_0} \right)^{u_{i}} x_i^{-u_{i+2} v_0}, & \\
  x_{i+1} \longmapsto  x_i^{v_0} x_{i+1}^{v_0 u_i} x_{i}^{-{v_0}^2 u_{i+1}}, &\\
    x_{i+2} \longmapsto  x_i^{v_0^2}. &\\
\end{array}
\right.
$$
Comparing the obtained automorphisms it can be seen that the long relation of the braid group holds when it acts on the generators of $F_n$. The fact that the relation holds with the action on subgroups generated by elements
$u_j$ and $v_j$ follows from the fact that the automorphisms $\varphi_M (\sigma_i)$ act as permutations.

We verify the correctness of the mixed relation
$$
\varphi_M  (\rho_i) \varphi_M (\rho_{i+1}) \varphi_M (\sigma_i) =
 \varphi_M (\sigma_{i+1}) \varphi_M (\rho_i) \varphi_M (\rho_{i+1}).
$$
For that purpose we find the  action of
$\varphi_M  (\rho_i) \varphi_M (\rho_{i+1}) \varphi_M (\sigma_i)$
on the generators $x_1$, $\ldots $, $x_n$. We have
$$
\varphi_M  (\rho_i) \varphi_M (\rho_{i+1}) \varphi_M (\sigma_i) :
\left\{
\begin{array}{ll}
  x_{i} \longmapsto  x_{i+2}^{v_i^{-1} v_{i+1}^{-1}}, & \\
  x_{i+1} \longmapsto  \left( x_i x_{i+1}^{u_i} x_{i}^{-u_{i+1} v_0} \right)^{v_{i+2}}, &\\
      x_{i+2} \longmapsto  x_i^{v_0 v_{i+2}}. &\\
\end{array}
\right.
$$
Analogously we find  the  action of
$\varphi_M (\sigma_{i+1}) \varphi_M (\rho_i) \varphi_M (\rho_{i+1})$. We obtain
$$
\varphi_M (\sigma_{i+1}) \varphi_M (\rho_i) \varphi_M (\rho_{i+1}) :
\left\{
\begin{array}{ll}
  x_{i} \longmapsto x_{i+2}^{v_{i+1}^{-1} v_i^{-1}}, & \\
  x_{i+1} \longmapsto  x_i^{v_{i+2}} x_{i+1}^{v_{i+2} u_i} x_{i}^{-v_{i+2} u_{i+1} v_0}, &\\
      x_{i+2} \longmapsto  x_i^{v_{i+2} v_0}. &\\
\end{array}
\right.
$$
Comparing the obtained automorphisms it is clear that the mixed relation of the braid group holds on the generators of $F_n$. The fact that the relation holds with the action on subgroups generated by elements
$u_j$ and $v_j$ follows because the automorphisms $\varphi_M (\sigma_i)$ and $\varphi_M (\rho_i)$ act as permutations.

The remaining relations are verified analogously. The theorem is proved.

\bigskip

It is easy to check that the element $\varphi_M (\sigma_i^{-1})$ acts on the generators of the group $F_n$ by the rule
$$
\varphi_M (\sigma_i^{-1}) :
\left\{
\begin{array}{ll}
  x_{i} \longmapsto {x_{i+1}}^{v_0^{-1}}, & \\
  x_{i+1} \longmapsto  \left( x_{i+1}^{-{v_0^{-1}}} x_{i} x_{i+1}^{u_{i}} \right)^{u_{i+1}^{-1}}. &\\
\end{array}
\right.
$$

We show that this representation generalizes the representations
 $\varphi_{SW}$ and $\varphi_{BD}$. It is true that

\medskip

{\bf Proposition 2.}
{\it 1) If we put
$$
v_1 = v_2 = \ldots = v_n = v
$$
in the group $F_{n,2n+1}$
and insert new generators $z_1,\ldots z_n$ such that
$$
z_1 = x_1, z_2 = x_2^v, z_3 = x_3^{v^2}, \ldots, z_n = x_n^{v^{n-1}},~~~w_i = v u_i,~~~i = 1, 2, \ldots, n,
$$
then we obtain the group $F_{n,n+1}$ on which the representation
$\varphi_M$ induces the representation $\varphi_{SW}$.

2) If we put
$$
u_1 = u_2 = \ldots = u_n = 1, v_0 = u,  v_1 = v_2 = \ldots = v_n = v
$$
in the group $F_{n,2n+1}$, then we obtain the group $F_{n,2}$ on which the representation
$\varphi_M$ induces the representation $\varphi_{BD}$. }

\medskip

{\bf Proof.} 1) Note that on the new generators $z_1,\ldots,z_n$ the automorphisms $\varphi_M(\sigma_i)$ and $\varphi_M (\rho_i)$ act by the rule
 $$
\varphi_M (\sigma_i) :
\left\{
\begin{array}{ll}
  z_{i} \longmapsto z_i z_{i+1}^{v u_i} z_i^{-u_{i+1}}, & \\
  z_{i+1} \longmapsto  z_{i}^{v^{-1}}, &\\
\end{array}
\right.~~~
\varphi_M (\rho_i) :
\left\{
\begin{array}{ll}
  z_{i} \longmapsto  z_{i+1}, & \\
  z_{i+1} \longmapsto  z_{i}. &\\
\end{array}
\right.
$$
Assuming $w_i = v u_i$ for $i = 1, 2, \ldots, n$, we obtain the Silver-Williams representation.

The second part of the proposition is verified analogously.
The proposition is proved.

\medskip

{\bf The representation of the virtual pure braid group.}
The virtual pure braid group $VP_n$ is the kernel of the homomorphism  $VB_n \longrightarrow S_n$ which acts on the generators:
$$
\sigma_i \longmapsto \rho_i,~~\rho_i \longmapsto \rho_i,~~i = 1, 2, \ldots, n-1.
$$
It was found out in \cite{B} that the group $VP_n$ is generated by elements
$$
\lambda_{i,i+1} = \rho_i \, \sigma_i^{-1},~~~
\lambda_{i+1,i} = \rho_i \, \lambda_{i,i+1} \, \rho_i = \sigma_i^{-1} \, \rho_i,
~~~i=1, 2, \ldots, n-1,
$$
$$
\lambda_{ij} = \rho_{j-1} \, \rho_{j-2} \ldots \rho_{i+1} \, \lambda_{i,i+1} \, \rho_{i+1}
\ldots \rho_{j-2} \, \rho_{j-1},
$$
$$
\lambda_{ji} = \rho_{j-1} \, \rho_{j-2} \ldots \rho_{i+1} \, \lambda_{i+1,i} \, \rho_{i+1}
\ldots \rho_{j-2} \, \rho_{j-1}, ~~~1 \leq i < j-1 \leq n-1.
$$
We find the images of the generators $\lambda_{ij}$, $\lambda_{ji}$ under the representation $\varphi_M$. It is true

\medskip

{\bf Proposition 3.}
{\it The representation $\varphi_M$ acts on the generators of  $VP_n$ by the formulae
$$
\varphi_M (\lambda_{i,i+1}) :
\left\{
\begin{array}{l}
  x_{i} \longmapsto \left( x_{i+1}^{-1} x_{i}  x_{i+1}^{u_{i+1}} \right)^{u_{i}^{-1} v_{i+1}^{-1}},  \\
  x_{i+1} \longmapsto   x_{i+1}^{v_i}, \\
\end{array}
\right.~~
\varphi_M (\lambda_{i+1,i}) :
\left\{
\begin{array}{l}
  x_{i} \longmapsto   x_{i}^{v_{i+1}}, \\
  x_{i+1} \longmapsto \left( x_{i}^{-v_{i+1}} x_{i+1}^{v_i^{-1}}  x_{i}^{v_{i+1} u_i} \right)^{u_{i+1}^{-1}},  \\
\end{array}
\right.
$$
$$
\varphi_M (\lambda_{ij}) :
\left\{
\begin{array}{l}
  x_{i} \longmapsto \left( x_{j}^{-v_{i+1}^{-1} \ldots v_{j-1}^{-1}} x_{i}   x_{j}^{-v_{i+1}^{-1} \ldots v_{j-1}^{-1} u_j} \right)^{u_{i}^{-1} v_{j}^{-1}},  \\
  x_{j} \longmapsto   x_{j}^{v_i}, \\
\end{array}
\right.
$$
$$
\varphi_M (\lambda_{ji}) :
\left\{
\begin{array}{l}
  x_{i} \longmapsto   x_{i}^{v_j}, \\
  x_{j} \longmapsto \left( x_{i}^{-v_j} x_{j}^{v_{i}^{-1} \ldots v_{j-1}^{-1}} x_{i}^{v_{j} u_i} \right)^{u_{j}^{-1} v_{i+1} \ldots v_{j-1}},  \\
\end{array}
\right.
$$
where $1 \leq i < j-1 \leq n-1$.}

\bigskip

{\bf The representation of the welded braid group.}
We consider $F_{n,n+1}$ which equals to the free product $ F_n * \mathbb{Z}^{n+1}$, where  $\mathbb{Z}^{n+1}$ is a free abelian group of the rank $n+1$ freely generated by elements $v, u_1, u_2, \ldots, u_n$. The following statement is true

\medskip

{\bf Proposition 4.}
{\it The mapping $\psi_M : WB_n \longrightarrow \mathrm{Aut}(F_{n,n+1})$ defined by the action on generators
$$
\psi_M (\sigma_i) :
\left\{
\begin{array}{ll}
  x_i \longmapsto  x_i x_{i+1}^{u_i} x_i^{- v u_{i+1} }, &\\
  x_{i+1} \longmapsto x_i^{v}, & \\
\end{array}
\right.~~~
\psi_M (\sigma_i) :
\left\{
\begin{array}{ll}
  u_i \longmapsto  u_{i+1}, &\\
  u_{i+1} \longmapsto u_i, & \\
\end{array}
\right.
$$
$$
\psi_M (\alpha_i) :
\left\{
\begin{array}{ll}
  x_i \longmapsto  x_{i+1}, &\\
  x_{i+1} \longmapsto x_i, & \\
\end{array}
\right.~~~
\psi_M (\alpha_i) :
\left\{
\begin{array}{ll}
  u_i \longmapsto  u_{i+1}, &\\
  u_{i+1} \longmapsto u_i & \\
\end{array}
\right.
$$
determines the representation of the group $WB_n$ into the group} $\mathrm{Aut}(F_{n,n+1})$.

It could be seen that the representation of the group $WB_n$ is obtained from the one of the group $VB_n$ providing all generators $v_i$ are sent to 1, where $1 \leq i \leq n$.

{\bf Proof.}
It is sufficient to check that the relation $\mathcal{F}_1$ holds
$$
\psi_M (\rho_i) \psi_M (\sigma_{i+1}) \psi_M (\sigma_i) =
 \psi_M (\sigma_{i+1}) \psi_M (\sigma_i) \psi_M (\rho_{i+1})
$$
which is the forbidden relation in $VB_n$.
We find the action of the automorphism $\psi_M (\rho_i) \psi_M (\sigma_{i+1}) \psi_M (\sigma_i)$ on the generators of  $F_n$. We obtain
$$
\psi_M (\rho_i) \psi_M (\sigma_{i+1}) \psi_M (\sigma_i) :
\left\{
\begin{array}{ll}
  x_{i} \longmapsto x_i  x_{i+2}^{u_i} x_{i}^{-u_{i+2}}, & \\
  x_{i+1} \longmapsto  x_i x_{i+1}^{u_i} x_{i}^{-u_{i+1}}, &\\
  x_{i+2} \longmapsto x_i. & \\
\end{array}
\right.
$$
Analogously we find the action of the automorphism
$\psi_M (\sigma_{i+1}) \psi_M (\sigma_i) \psi_M (\rho_{i+1})$. We have
$$
\psi_M (\sigma_{i+1}) \psi_M (\sigma_i) \psi_M (\rho_{i+1}) :
\left\{
\begin{array}{ll}
  x_{i} \longmapsto x_i  x_{i+2}^{u_i} x_{i}^{-u_{i+2}}, & \\
  x_{i+1} \longmapsto  x_i x_{i+1}^{u_i} x_{i}^{-u_{i+1}}, &\\
   x_{i+2} \longmapsto x_i. & \\
\end{array}
\right.
$$
Comparing the obtained automorphisms we see  that the required relation holds.

The proposition is proved.

\bigskip

As it is known the relation $\mathcal{F}_1$ holds in $WB_n$ but at the same time the relation $\mathcal{F}_2$ does not hold.  It would be interesting to construct  representations of other factor-groups of the group $VB_n$. For example, the representation of the group $UVB_n$ in which both forbidden relations hold.
We verify the conditions for the second forbidden relation $\mathcal{F}_2$ to hold:
$$
\psi_M (\rho_{i+1}) \psi_M (\sigma_{i}) \psi_M (\sigma_{i+1}) =
 \psi_M (\sigma_{i}) \psi_M (\sigma_{i+1}) \psi_M (\rho_{i}).
$$
We find the action of the automorphism $\psi_M (\rho_{i+1}) \psi_M (\sigma_{i}) \psi_M (\sigma_{i+1})$
on the generators of  $F_n$. We obtain
$$
\psi_M (\rho_{i+1}) \psi_M (\sigma_{i}) \psi_M (\sigma_{i+1}) :
\left\{
\begin{array}{ll}
  x_{i} \longmapsto x_i \left( x_{i+1} x_{i+2}^{u_{i+1}} x_{i+1}^{-u_{i+2}} \right)^{u_i} x_{i}^{-u_{i+2}}, & \\
  x_{i+1} \longmapsto   x_{i+1}, &\\
    x_{i+2} \longmapsto   x_{i}. &\\
\end{array}
\right.
$$
Analogously we find the action of the automorphism $\psi_M (\sigma_{i}) \psi_M (\sigma_{i+1}) \psi_M (\rho_{i})$.
We have
$$
\psi_M (\sigma_{i}) \psi_M (\sigma_{i+1}) \psi_M (\rho_{i}) :
\left\{
\begin{array}{ll}
  x_{i} \longmapsto x_{i+1}  \left( x_{i} x_{i+2}^{u_{i}} x_{i}^{-u_{i+2}} \right)^{u_{i+1}} x_{i+1}^{-u_{i+2}}, & \\
  x_{i+1} \longmapsto   x_{i+1}, &\\
    x_{i+2} \longmapsto   x_{i}. &\\
\end{array}
\right.
$$
Comparing the obtained automorphisms it can be seen that the required relation holds providing the following restriction
$$
\left( x_{i+1} x_{i+2}^{u_{i+1}} x_{i+1}^{-u_{i+2}} \right)^{u_i} x_{i}^{-u_{i+2}} =    x_{i+1}  \left( x_{i} x_{i+2}^{u_{i}} x_{i}^{-u_{i+2}} \right)^{u_{i+1}} x_{i+1}^{-u_{i+2}},
$$
i.~e.
$$
\left( x_{i+1}^{-u_i} x_{i}^{-1} x_i x_{i}^{u_{i+1}} \right) \cdot x_{i+2}^{ u_i u_{i+1}} =
 x_{i+2}^{u_{i} u_{i+1}} \cdot \left( x_{i+1}^{-u_i} x_{i}^{-1} x_i x_{i}^{u_{i+1}} \right)^{u_{i+2}}.
$$
These are the relations that are of the form $a b = b a^c$ or equally $[ac, b] = 1$.

\medskip

{\bf Question 1.} What is the group where these relations hold?

\newpage


\begin{center}
{\bf \S~3. Representation that is equivalent to the representation $\varphi_M$}
\end{center}

\vspace{0.5cm}

The constructed representation $\varphi_M$ is not an extension of the Artin representation which is known to be faithful on the subgroup $B_n$.
We show that the representation $\varphi_M$ is equivalent to the simpler one which is an extension of the Artin representation.
Let $F_{n,n} = F_n * \mathbb{Z}^n$ where $F_n = \langle y_1, y_2, \ldots, y_n \rangle$ is the free group  and $\mathbb{Z}^n = \langle v_1, v_2, \ldots, v_n \rangle$ is the free abelian group of the rank $n$.
It is true

\medskip

{\bf Theorem 2.} {\it The representation $\widetilde{\varphi}_M : VB_n \longrightarrow \mathrm{Aut}(F_{n,n})$ defined by the action on the generators
$$
\widetilde{\varphi}_{M}(\sigma_i) :
\left\{
\begin{array}{l}
  y_i \longmapsto  y_i y_{i+1} y_i^{-1}, \\
  y_{i+1} \longmapsto y_i,  \\
\end{array}
\right.~~~
\widetilde{\varphi}_{M}(\sigma_i) :
\left\{
\begin{array}{l}
  v_i \longmapsto  v_{i+1}, \\
  v_{i+1} \longmapsto v_i,  \\
\end{array}
\right.
$$
$$
\widetilde{\varphi}_{M}(\rho_i) :
\left\{
\begin{array}{l}
  y_i \longmapsto  y_{i+1}^{v_i^{-1}}, \\
  y_{i+1} \longmapsto y_i^{v_{i+1}},  \\
\end{array}
\right.~~~
\widetilde{\varphi}_{M}(\rho_i) :
\left\{
\begin{array}{l}
  v_i \longmapsto  v_{i+1}, \\
  v_{i+1} \longmapsto v_i  \\
\end{array}
\right.
$$
is equivalent to the representation $\varphi_M$.}

\medskip

The theorem proof is due to the following lemmas.

\medskip

{\bf Lemma 1. }
{\it If we take new generators $y_i$ and $w_i$ in the group $F_{n,2n+1}$ such that
$y_i=(x_i u_i^{-1}v_0^{-1})^{v_0^{-(i-1)}}$, $w_i=v_iv_0^{-1}$, $i=1,\ldots,n$,
instead of generators $x_i $, $v_i$,  $i=1,\ldots,n$, than the representation
$$
\varphi_{M} : VB_n \longrightarrow \mathrm{Aut}(F_{n,2n+1})
$$
acts on new generators
$y_1,\ldots,y_n$, $u_1,\ldots,u_n$, $v_0,w_1,\ldots,w_n$ by the following way
$$
\varphi_{M}(\sigma_i) :
\left\{
\begin{array}{l}
  y_i \longmapsto  y_i y_{i+1} y_i^{-1}, \\
  y_{i+1} \longmapsto y_i,  \\
\end{array}
\right.~~~
\varphi_{M}(\sigma_i) :
\left\{
\begin{array}{l}
  u_i \longmapsto  u_{i+1}, \\
  u_{i+1} \longmapsto u_i,  \\
\end{array}
\right.
~~~
\varphi_{M}(\sigma_i) :
\left\{
\begin{array}{l}
  w_i \longmapsto  w_{i+1}, \\
  w_{i+1} \longmapsto w_i,  \\
\end{array}
\right.
$$
$$
\varphi_{M}(\rho_i) :
\left\{
\begin{array}{l}
  y_i \longmapsto  y_{i+1}^{w_i^{-1}}, \\
  y_{i+1} \longmapsto y_i^{w_{i+1}},  \\
\end{array}
\right.~~~
\varphi_{M}(\rho_i) :
\left\{
\begin{array}{l}
  u_i \longmapsto  u_{i+1}, \\
  u_{i+1} \longmapsto u_i,  \\
\end{array}
\right.
~~~
\varphi_{M}(\rho_i) :
\left\{
\begin{array}{l}
  w_i \longmapsto  w_{i+1}, \\
  w_{i+1} \longmapsto w_i.  \\
\end{array}
\right.
$$
}

{\bf Proof.}
It is sufficient to find the action $\varphi_{M}(\sigma_i)$ and $\varphi_{M}(\rho_i)$
on $y_i$, $y_{i+1}$, $w_i$, $w_{i+1}$, $i=1,\ldots,n-1$. We obtain
$$
\begin{array}{rl}
  \varphi_{M}(\sigma_i)(y_i) & =  (x_i u_i^{-1}v_0^{-1})^{v_0^{-(i-1)}}
                        \longmapsto (x_i x_{i+1}^{u_i}x_i^{-v_0u_{i+1}}u_{i+1}^{-1}v_0^{-1} )^{v_0^{-(i-1)}}\\
       &=(x_iu_i^{-1} x_{i+1}u_iv_0^{-1} u_{i+1}^{-1} x_i^{-1}v_0u_{i+1}u_{i+1}^{-1}v_0^{-1} )^{v_0^{-(i-1)}}\\
       &=(x_iu_i^{-1} x_{i+1}u_{i+1}^{-1}u_iv_0^{-1} x_i^{-1})^{v_0^{-(i-1)}}\\
       &=(x_iu_i^{-1}v_0^{-1} v_0x_{i+1}u_{i+1}^{-1}v_0^{-1}v_0^{-1}v_0u_i x_i^{-1})^{v_0^{-(i-1)}}\\
       &=(x_iu_i^{-1}v_0^{-1})^{v_0^{-(i-1)}}
              (x_{i+1}u_{i+1}^{-1})^{v_0^{-i}}
                    (x_iu_i^{-1}v_0^{-1})^{-v_0^{-(i-1)}} \\
       &= y_iy_{i+1} y_i^{-1}, \\
  \varphi_{M}(\sigma_i)(y_{i+1}) & =  (x_{i+1} u_{i+1}^{-1}v_0^{-1})^{v_0^{-i}}
                        \mapsto (x_i^{v_0} u_i^{-1}v_0^{-1} )^{v_0^{-i}}=
                        (x_iu_i^{-1} v_0^{-1} )^{v_0^{-(i-1)}}= y_i, \\
  \varphi_{M}(\sigma_i)(w_i) & =  v_iv_0^{-1} \mapsto v_{i+1}v_0^{-1}=w_{i+1}, \\
  \varphi_{M}(\sigma_i)(w_{i+1}) & =  v_{i+1}v_0^{-1} \mapsto v_iv_0^{-1}=w_i, \\
\end{array}
$$
$$
\begin{array}{ll}
 \varphi_{M}(\rho_i)(y_i)      & =  (x_i u_i^{-1}v_0^{-1})^{v_0^{-(i-1)}}
                        \longmapsto (x_{i+1}^{v_i^{-1}} u_{i+1}^{-1}v_0^{-1} )^{v_0^{-(i-1)}}\\
                              &=(x_{i+1} u_{i+1}^{-1}v_0^{-1} )^{v_0^{-i} v_i^{-1} v_0 }=y_{i+1}^{w_i^{-1}},\\
  \varphi_{M}(\rho_i)(w_i)     & =  v_iv_0^{-1} \mapsto v_{i+1}v_0^{-1}=w_{i+1}, \\
  \varphi_{M}(\rho_i)(w_{i+1}) & =  v_{i+1}v_0^{-1} \mapsto v_iv_0^{-1}=w_i. \\
\end{array}
$$

The lemma is proved.

\medskip

On account of the Lemma 1 the subgroups
$$
F_{n,n}=\left\langle y_1,\ldots,y_n, w_1,\ldots,w_n \right\rangle, \quad
W=\left\langle  w_1,\ldots,w_n \right\rangle, \quad
U=\left\langle  u_1,\ldots,u_n \right\rangle, \quad
V_0=\left\langle  v_0 \right\rangle
$$
of  $F_{n,2n+1}$ are invariants under the action of  $\varphi_M (VB_n)$. In particular, the representation $\varphi_M$ induces the representation $\widetilde{\varphi}_M$ of the group $VB_n$ into the group $Aut (F_{n,n})$.
Also note that the restrictions $\varphi_M(VB_n)$ on the subgroups $W$ and $U$ are equivalent with respect to the isomorphisms $u_i\leftrightarrows w_i$ where
$i=1,\ldots,n$.

\medskip

{\bf Lemma 2. }
{\it The representations $\varphi_M$ and $\widetilde{\varphi}_M$ are equivalent, i. ~e. $\varphi_M(\beta) = 1$ if and only if $\widetilde{\varphi}_M(\beta) = 1$ for some $\beta \in VB_n$.}

{\bf Proof} follows directly from the fact that the restrictions of $\varphi_M(VB_n)$ on the subgroups $W$ and $U$ are equivalent. The lemma is proved.

\medskip

Theorem 2 follows from the proved lemmas.

\medskip

Hence the most general representation $\widetilde{\varphi}_M$ is in fact the extension of the Artin representation of the braid group $B_n$.

Using the result of O. Chterental \cite{Ch} who proved that the representation $\varphi_A$ of the group $VB_n$ is not faithful for $n \geq 4$, it could be shown that the representations $\varphi_{SW}$ and $\varphi_{BD}$ are not faithful for $n \geq 4$ as well.
It is true

\medskip

{\bf Proposition 5.} {\it The element $\beta = (\sigma_2^{-1} \rho_1 \sigma_2 \rho_3)^3 \in VP_4$ lies in the kernel of the representations $\varphi_{SW}$ as well as $\varphi_{BD}$ and it does not lie in the kernel of the representation} $\widetilde{\varphi}_M$.

\medskip

{\bf Proof.} The straightforward verification that $\varphi_{SW}(\beta) = 1$ and $\varphi_{BD}(\beta) = 1$ gives us the first part of the proof.

Further we denote $ \sigma_2^{-1} \rho_1 \sigma_2 \rho_3$ by $\alpha$, then $\beta = \alpha^3$.
To check the remaining part of the proof we find the automorphism $\widetilde{\varphi}_M(\alpha)$. We obtain
$$
\widetilde{\varphi}_{M}(\alpha) :
\left\{
\begin{array}{l}
  x_1 \longmapsto  (x_2 x_{4}^{v_3^{-1}} x_2^{-1})^{v_1^{-1}}, \\
  x_2 \longmapsto x_2,  \\
    x_3 \longmapsto  x_2^{-1} x_{1}^{v_4} x_2, \\
  x_4 \longmapsto x_3^{v_4},~~ \\
\end{array}
\right.~~~
\widetilde{\varphi}_{M}(\alpha) :
\left\{
\begin{array}{l}
  u_1 \longmapsto  u_{4}, \\
  u_{2} \longmapsto u_2,  \\
  u_3 \longmapsto u_1, \\
    u_4 \longmapsto u_3. \\
\end{array}
\right.
$$
Next we find
$$
\widetilde{\varphi}_{M}(\beta) :
\left\{
\begin{array}{l}
  x_1 \longmapsto  x_2^{v_3^{-1}} \left( x_2^{-1} x_{1}^{v_4} x_2 \right)^{v_4^{-1}} x_2^{-v_3^{-1}}, \\
  x_2 \longmapsto x_2,  \\
    x_3 \longmapsto  x_2^{-1} x_{2}^{v_1} x_3 x_2^{-v_4^{-1}} x_2, \\
  x_4 \longmapsto x_2^{-v_1} \left( x_2 x_{4}^{v_3^{-1}} x_2^{-1} \right)^{v_3} x_2^{v_1},~~ \\
\end{array}
\right.
$$
i.~e. $\widetilde{\varphi}_{M}(\beta) \not= 1$. The proposition is proved.

\medskip

Nevertheless the following question remains open

{\bf Question 2.} Is the representation $\varphi_{A} : VB_n \longrightarrow \mathrm{Aut}(F_{n+1})$ faithful for $n = 3$?

\vspace{0.8cm}

\begin{center}
{\bf \S~4. The link groups}
\end{center}

\vspace{0.5cm}

To find a presentation of the link group two approaches could be used. The first one is the Wirtinger approach when the link diagram on the plane  is used to find the generators and defining relations of the link group. The second one is to consider the link as a closed braid and  use the braid group representation into the automorphism group of some group. Both these approaches are described in the book \cite{Bir} for classical links and in the paper \cite{BB} for virtual ones. In particular, for each virtual braid $\beta$ the groups $G_A(\beta)$,  $G_{SW}(\beta)$, $G_{BD}(\beta)$ and $G_M(\beta)$ could be defined using the representations $\varphi_A$, $\varphi_{SW}$, $\varphi_{BD}$ and $\varphi_M$ correspondingly. In addition it should be proved that the corresponding group depends only on the link $\widehat{\beta}$ (the link $\widehat{\beta}$ is the closure of the braid $\beta$) if we want to construct the link invariant.

\begin{figure}[h]
\noindent\centering{\includegraphics[height=0.3 \textwidth]{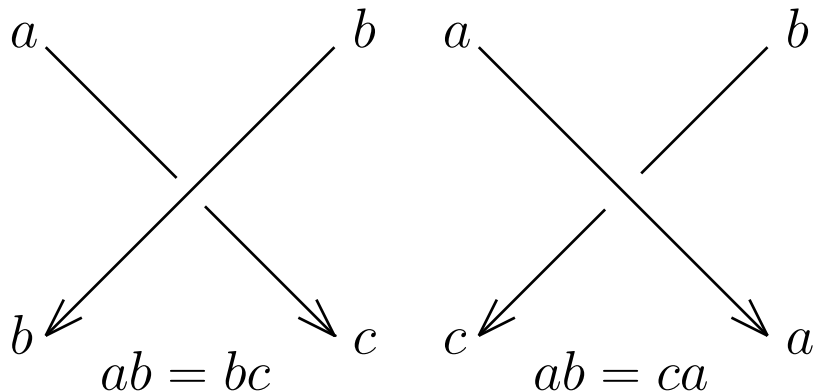}}
\caption{Relations in classical crossings}
\label{sigma}
\end{figure}

{\bf The classical link group.}
Recall that for a classical link $L$ in $\mathbb{S}^3$ \textit{its group} $G(L)$ is the fundamental group $\pi_1(\mathbb{S}^3 \setminus L)$
of the complement of $L$ in $\mathbb{S}^3$.
To find  generators and relations of the group $G(L)$ either the link diagram or the consideration of the link as the closed braid and the Artin representation could be used. In the case we use the link diagram the generator is assigned to each connected component and the following relation is assigned to each classical crossing:
$a b = b c$ or its equivalent  $c = b^{-1} a b$ for the positive crossing
(see Fig. \ref{sigma} on the left) and $a b = c a$ or its equivalent $c =  a b a^{-1}$ for the negative crossing (see Fig. \ref{sigma} on the right).

{\bf The generalized Alexander group.}
D. Silver and S. Williams \cite{SW} defined the generalized Alexander  group $\widetilde{A}_L$ for the virtual $d$-component link
$L = K_1 \cup K_2 \cup \ldots \cup K_d$.  For that purpose the free abelian group $\Pi_d = \langle u_1, u_2, \ldots, u_d, v \rangle$ of the rank $d+1$ was used.
In particular provided $L$ a knot, the free abelian group $\Pi_1 = \langle u, v \rangle$ of the rank $2$ is used.

\begin{figure}[h]
\noindent\centering{\includegraphics[height=0.3 \textwidth]{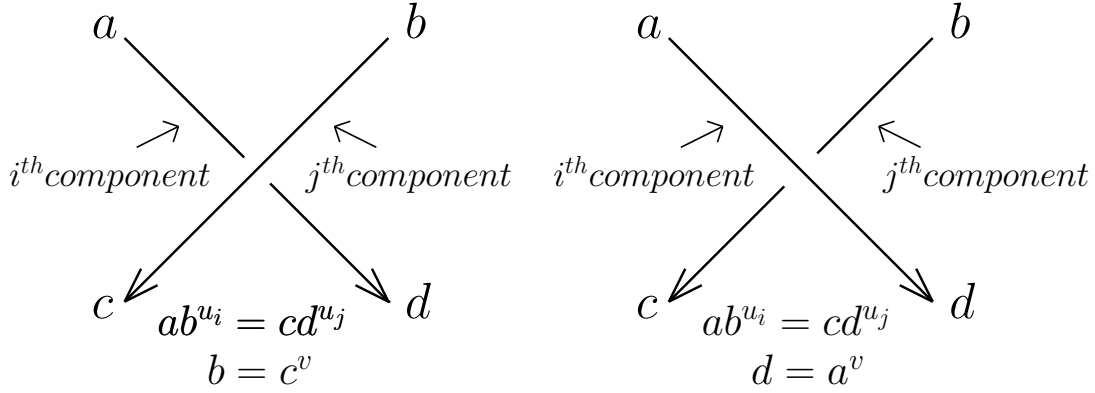}}
\caption{Relation in the Alexander generalized group}
\label{sigmaSW}
\end{figure}
The particular generator is assigned to each arc of the diagram $L$ from one classical crossing to another one.
Two relations are assigned to each positive crossing:
$a b^{u_i} = c d^{u_j}$, $b = c^v$ (see Fig. \ref{sigmaSW} on the left) which could be written in the form
$$
c = b^{v^{-1}},~~d = b^{-v^{-1} u_j^{-1}} a^{u_j^{-1}} b^{u_i u_j^{-1}}.
$$
Two relations are assigned to each negative crossing:
$a b^{u_i} = c d^{u_j}$, $d = a^v$ (see Fig. \ref{sigmaSW} on the right) which could be written in the form
$$
c = a b^{u_i} a^{-v u_j},~~d =  a^{v}.
$$
No relation is assigned to the virtual crossing.

{\bf The welded braid group.} For the welded $d$ component link $L$ the generalized Alexander group $A_L$ is constructed (see \cite{SW}) which is obtained from the corresponding virtual link by putting $v = 1$,
i. e. the following relations correspond to each classical crossing
$$
a b^{u_i} = c d^{u_j},~~b = c
$$
for the positive case and the relations
$$
a b^{u_i} = c d^{u_j},~~d = a
$$
for the negative one.

The generalized Alexander group of the trivial $d$ component link
$U_d$ is a free product $F_d * \mathbb{Z}^d$ of the free group of the rank $d$ and the free abelian group of the same rank.

\newpage

\begin{center}
{\bf \S~5. The virtual link groups}
\end{center}

\vspace{0.5cm}

In the present section the new definition of the virtual link group is given.
For this two approaches are used: the group will be defined at first with the help of the representation $\varphi_M$ (the braid approach) then through the virtual diagram. Next the equivalence of these two definitions is shown.

{\bf The braid approach.} We describe the general way which allows to construct the link invariants through the representation of the virtual braid group by automorphisms of some group. Assume that we have a representation $\varphi : VB_n \longrightarrow \mathrm{Aut}(H)$
of the virtual braid group into the automorphism group of some group
$H = \langle h_1, h_2, \ldots, h_m ~\|~ \mathcal{R}\rangle$, where $\mathcal{R}$ is the set of defining relations. The following group is assigned to the virtual braid $\beta \in VB_n$:
$$
G_{\varphi}(\beta) = \langle h_1, h_2, \ldots, h_m~\|~
\mathcal{R}, h_i = \varphi (\beta) (h_i),  ~~i = 1, 2, \ldots, m \rangle.
$$
The group $G_{\varphi}$ is the invariant of the virtual link providing the verification that the group $G_{\varphi}(\beta)$ is isomorphic to $G_{\varphi}(\beta')$ for each braid $\beta'$ such that the links $\widehat{\beta}$ and $\widehat{\beta'}$ are equivalent.

This approach is used for the previously defined representation $\varphi_M$.
Given $\beta \in VB_n$, the \textit{ group of the braid} $\beta$ is  the following group
\begin{align*}
\begin{split}
G_M(\beta) &=
 \langle x_1, x_2,\ldots, x_n, u_1, u_2,\ldots, u_n, v_0, v_1,  \ldots, v_n ~|| ~
 [u_i,u_j]=[v_k,v_l]=[u_i,v_k]=1, \\
x_i = \varphi_{M} (\beta)& (x_i),  ~~u_i = \varphi_{M} (\beta) (u_i),
~~v_i = \varphi_{M} (\beta) (v_i), ~~i,j = 1, 2, \ldots, n , ~~k,l = 0, 1, \ldots, n\rangle.
\end{split}
\end{align*}


Let show that the group defined in this manner is indeed the link invariant. More precisely it is true

\medskip

{\bf Theorem 3.}
{\it Given $\beta \in VB_n$ and $\beta' \in VB_m$ the two virtual braids such that theirs closures define the same link $L$,
 then $G_M(\beta) \cong G_M(\beta')$,
i.~e. the group $G_M(\beta)$ is the link $L$ invariant and it could be denoted by the symbol $G_M(\widehat{\beta}) = G_M(L)$.}

\medskip

The next paragraph is dedicated to the proof of the theorem.

\bigskip

{\bf The diagram approach.} Now we show the way of how to find the virtual link group using its diagram. Let $D_L$ be a diagram of a virtual $d$-component link $L$. This diagram is divided into arcs from one crossing (classical or virtual) to other and  for every obtained arc a generator is assigned. Thus we have some set of generators $a, b, c, \ldots$.
We numerate all components of $D_L$ by numbers from 1 to $d$ and assign to the $i$ th component two generators: $u_i$ and $v_i$.

\begin{figure}[h]
\noindent\centering{\includegraphics[height=0.3 \textwidth]{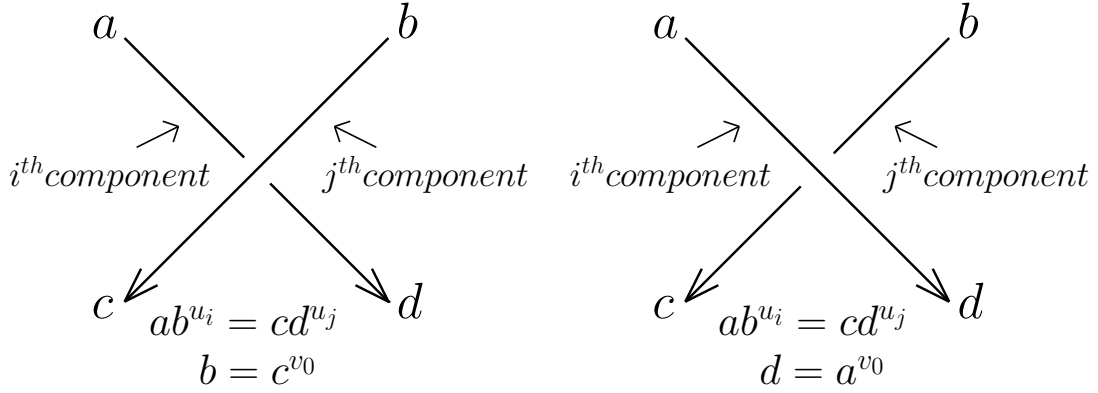}}
\caption{The relations in classical crossings of the group $G_M$}
\label{sigmaM}
\end{figure}

\begin{figure}[h]
\noindent\centering{\includegraphics[height=0.35 \textwidth]{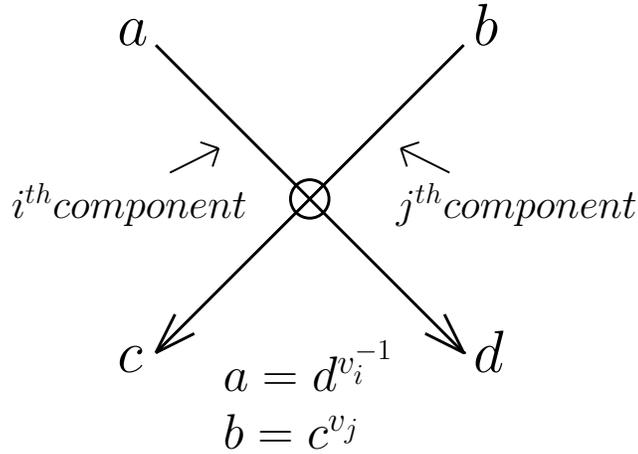}}
\caption{The relations in the virtual crossing of the group $G_M$}
\label{roM}
\end{figure}

A {\it group of the diagram} $D_L$ is the group generated by elements $a, b, c, \ldots$, $u_i$, $v_0$, $v_i$, $i = 1, \ldots, d$, and defined by the following system of relations: two relations correspond to each crossing(see Fig. \ref{sigmaM}, \ref{roM}).
Furthermore, all elements $u_i$ and $v_j$ commute pairwise , i. e. $\langle u_1, \ldots, u_d, v_0, v_1, \ldots, v_d \rangle$ is a free abelian group of the rank $2d+1$.

This group is denoted by $G_D(D_L)$. If $D_L'$ is another diagram of $L$, then $D_L$ can be transformed to $D_L'$ with the help of a finite number of generalized Reidemeister moves and a flat isotopy. Analogously to the case of classical link diagrams it is not difficult to prove

\medskip

{Proposition 6.} {\it If $D_L$ and $D_L'$ are two diagrams of a virtual link $L$, then groups $G_D(D_L)$ and $G_D(D_L')$ are isomorphic. Hence $G_D(D_L)$ is an invariant of $L$ and we denote it by $G_D(L)$.}

\medskip

Now we show that the group of a virtual link $L$ constructed by the braid is isomorphic to the group constructed by the diagram. More precisely, it is true

\medskip

{Proposition 7.} {\it Let $L$ be a virtual link, $D_L$ be its diagram, $\beta$ be a braid such that its closure $\widehat{\beta}$ is equivalent to $L$,  then $G_D(D_L)$ is isomorphic to $G_M(\beta)$.}

\medskip

{Proof.} Let $\beta$ be a virtual braid from the group $VB_n$ for some $n$. We consider the braid $\beta$ as a word in the alphabet $\{ \sigma_1^{\pm 1}, \sigma_2^{\pm 1}, \ldots, \sigma_{n-1}^{\pm 1}, \rho_1, \rho_2, \ldots, \rho_{n-1} \}$. Every generator $\sigma_i, \sigma_i^{-1}$ and $\rho_i$ has a geometric interpretation and can be presented by its diagram. Hence the braid
 $\beta$ can be presented by a diagram which is divided  in the layers so that each layer contains only one crossing: classical or virtual.
 The generators $\sigma_i$ or $\sigma_i^{-1}$ correspond to the classical crossing and the generator $\rho_i$ corresponds to the virtual one.
If the braid $\beta$ is  a  word of the length $m$  in the group $VB_n$, then its diagram contains $m$ layers. We consider the first layer. It contains $n$ arcs. The middle line of the first layer divides these arcs in half. We assign to starting parts of these arcs the symbols $x_{1,j}, u_{1,j}, v_{1,j}$, $j = 1, 2, \ldots, n$. To the next parts of the arcs which end  at the middle of the second layer we assign the symbols $x_{2,j}, u_{2,j}, v_{2,j}$, $j = 1, 2, \ldots, n$. By continuing this process  we divide the diagram $\beta$ into $n (m+1)$ arcs.

We are moving alongside the braid diagram from the bottom upwards, i.e. we are reading the word $\beta$ from right to left.
Let  the  $k$-th layer contain  the classical crossing, which corresponds to the elementary braid $\sigma_i$. In this layer we have two sets of arcs with the corresponding sets  of symbols $x_{kj}$, $u_{kj}$, $v_{kj}$
and $x_{k+1,j}$, $u_{k+1,j}$, $v_{k+1,j}$, $j=1,\ldots ,n$.
\begin{figure}[h]
\noindent\centering{\includegraphics[height=0.23 \textwidth]{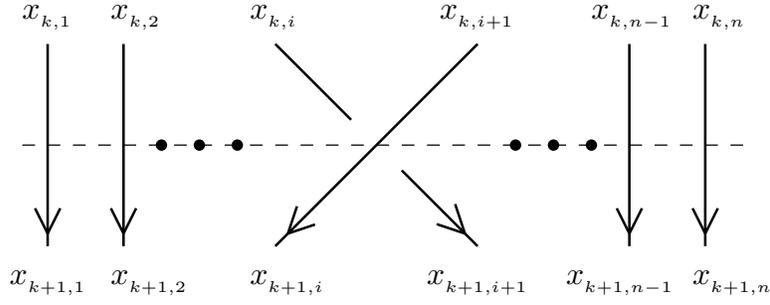}}
\caption{Classical positive crossing}
\label{sigmaDiag}
\end{figure}

Due to the formulae of the Theorem 1, the following automorphism corresponds to the
 braid $\sigma_i$
$$
\varphi_{M}(\sigma_i) :
\left\{
\begin{array}{l}
  x_{k+1,i} \longmapsto  x_{k,i} x_{k,i+1}^{u_{ki}} x_{k,i}^{-v_0u_{k,i+1}}, \\
  x_{k+1,i+1} \longmapsto x_{k,i}^{v_0}, \\
  x_{k+1,j} \longmapsto x_{k,j}, \\
\end{array}
\right.
\left\{
\begin{array}{l}
  u_{k+1,i}   \longmapsto  u_{k,i+1}, \\
  u_{k+1,i+1} \longmapsto  u_{k,i}, \\
  u_{k+1,j} \longmapsto  u_{k,j}, \\
\end{array}
\right.
\left\{
\begin{array}{l}
  v_{k+1,i}   \longmapsto  v_{k,i+1}, \\
  v_{k+1,i+1} \longmapsto  v_{k,i}, \\
  v_{k+1,j} \longmapsto  v_{k,j}, \\
\end{array}
\right.
$$
where  $j\neq i, i+1$. The case when the $k$-th layer contains $\sigma_i^{-1}$ is considered by analogy.

If the $k$-th layer contains  $\rho_i$ (see Fig. 6), then we have the following automorphism
\begin{figure}[h]
\noindent\centering{\includegraphics[height=0.23 \textwidth]{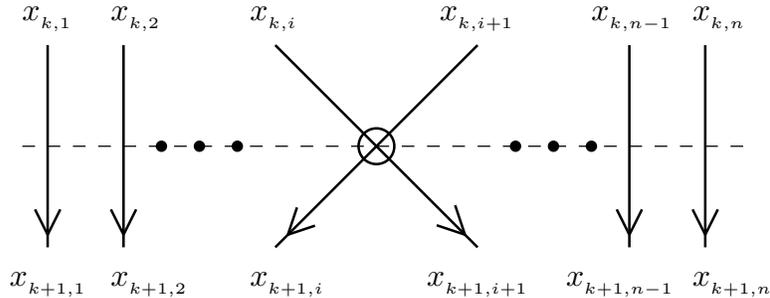}}
\caption{Virtual crossing}
\label{roDiag}
\end{figure}

$$
\varphi_{M}(\rho_i) :
\left\{
\begin{array}{l}
  x_{k+1,i} \longmapsto  x_{k,i+1}^{v_{ki}^{-1}}, \\
  x_{k+1,i+1} \longmapsto x_{k,i}^{v_{k,i+1}},  \\
  x_{k+1,j} \longmapsto x_{k,j},  \\
\end{array}
\right.
\left\{
\begin{array}{l}
  u_{k+1,i}   \longmapsto  u_{k,i+1}, \\
  u_{k+1,i+1} \longmapsto  u_{k,i}, \\
  u_{k+1,j} \longmapsto  u_{k,j}, \\
\end{array}
\right.
\left\{
\begin{array}{l}
  v_{k+1,i}   \longmapsto  v_{k,i+1}, \\
  v_{k+1,i+1} \longmapsto  v_{k,i}, \\
  v_{k+1,j} \longmapsto  v_{k,j}, \\
\end{array}
\right.
$$
where $j\neq i, i+1$.

Since in every layer top and bottom arcs are identified we can write the set of relations

1)
$$
\left\{
\begin{array}{l}
  x_{k+1,i}   =  x_{k,i} x_{k,i+1}^{u_{ki}} x_{k,i}^{-v_0u_{k,i+1}}, \\
  x_{k+1,i+1} = x_{k,i}^{v_0}, \\
  x_{k+1,j}   = x_{k,j},  \\
\end{array}
\right.
\left\{
\begin{array}{l}
  u_{k+1,i}   =  u_{k,i+1}, \\
  u_{k+1,i+1} =  u_{k,i}, \\
  u_{k+1,j} =  u_{k,j}, \\
\end{array}
\right.
\left\{
\begin{array}{l}
  v_{k+1,i}   =  v_{k,i+1}, \\
  v_{k+1,i+1} =  v_{k,i}, \\
  v_{k+1,j}   =  v_{k,j}, \\
\end{array}
\right.
$$
where $j\neq i, i+1$,
for  the $k$-th layer corresponded to $\sigma_i$;

2)
$$
\left\{
\begin{array}{l}
  x_{k+1,i}   =  x_{k,i+1}^{v_{ki}^{-1}}, \\
  x_{k+1,i+1} = x_{k,i}^{v_{k,i+1}},  \\
  x_{k+1,j} = x_{k,j},  \\
\end{array}
\right.
\left\{
\begin{array}{l}
  u_{k+1,i}   =  u_{k,i+1}, \\
  u_{k+1,i+1} =  u_{k,i}, \\
  u_{k+1,j} =  u_{k,j}, \\
\end{array}
\right.
\left\{
\begin{array}{l}
  v_{k+1,i}   =  v_{k,i+1}, \\
  v_{k+1,i+1} =  v_{k,i}, \\
  v_{k+1,j} =  v_{k,j}, \\
\end{array}
\right.
$$
where $j\neq i, i+1$,
for $k$-th layer  which corresponds to $\rho_i$.

Also we have to add

3) commutation relations:
$$
[u_{kl},u_{pq}]=[u_{kl},v_{pq}]=[v_{kl},v_{pq}]=[u_{kl},v_0]=[v_{kl},v_0]=1, ~~~k,p = 1, \ldots ,m+1, l, q = 1, \ldots, n.
$$

Using these relations, we move alongside the braid $\beta$ from the bottom upwards and we express the generators $x_{m+1,j}$, $u_{m+1,j}$, $v_{m+1,j}$ from the $m$-th layer through  the generators $x_{m,j}$, $u_{m,j}$, $v_{m,j}$, $j = 1, 2, \ldots, n$.
Next we consider the $m-1$-th layer and expressing the generators $x_{m+1,j}$, $u_{m+1,j}$, $v_{m+1,j}$ through  the generators
$x_{m-1,j}$, $u_{m-1,j}$, $v_{m-1,j}$, $j = 1, 2, \ldots, n$, and so on. And we express the generators $x_{m+1,j}$, $u_{m+1,j}$, $v_{m+1,j}$, through  the generators $x_{1,j}$, $u_{1,j}$, $v_{1,j}$, $j = 1, 2, \ldots, n$ when we arrive to the first layer.

To come to the
 closed braid $\widehat{\beta}$ we must connect initial  arcs of the diagram $\beta$ with the terminal arcs of this diagram, i.e. we must add the relations

 4)
$$
x_{m+1,j}=x_{1,j}, \quad u_{m+1,j}=u_{1,j}, \quad v_{m+1,j}=v_{1,j}, \quad j=1,\ldots,n.
$$

It follows from relations 1), 2), 4) that the generators $u_{k,l}$, $v_{k,l}$ that correspond to one component of the link are equal. This corresponds to the fact that under the homomorphism of $VB_n$ onto $S_n$ the image $\beta$ in $S_n$ is a product of $d$ independent cycles. Hence among the generators $u_{i,l}$ (also as among the generators $v_{i,l}$) there are only $d$ different ones.

Excluding transitional letters $x_{k,l}$, $u_{k,l}$, $v_{k,l}$, $k\geq 2$, we obtain relations
$$
\varphi_{M}(\beta)(x_1)=x_1, \ldots , \varphi_{M}(\beta)(x_n)=x_n,
$$
$$
\varphi_{M}(\beta)(u_1)=u_1, \ldots , \varphi_{M}(\beta)(u_n)=u_n,
$$
$$
\varphi_{M}(\beta)(v_1)=v_1, \ldots , \varphi_{M}(\beta)(v_n)=v_n,
$$
of the virtual link group constructed using the braid $\beta$.
Here  $x_{1,1}=x_1$, $\ldots $, $x_{1,n}=x_n$, $u_{1,1}=u_1$, $\ldots $, $u_{1,n}=u_n$, $v_{1,1}=v_1$, $\ldots $, $v_{1,n}=v_n$. Hence we get the representation of the group $G_M(\beta)$. The proposition is proved.

\medskip

{\bf Welded link group.}
A welded ling group could be defined analogously using the representation $\psi_M$.
Let a welded link $L_w$ be equivalent to a closure of the braid $\beta_w \in WB_n$
for some $n$. \textit{The group of the braid $\beta_{w}$} is  the following group
$$
G(\beta_w) = \langle x_1, x_2, \ldots, x_n, u_1, u_2, \ldots, u_n~||~ [u_i,u_j]=1,
$$
$$
x_i = \psi_M (\beta_w) (x_i),~~ u_i = \psi_M (\beta_w) (u_i), ~~i,j = 1, 2, \ldots, n \rangle.
$$
Taking into account the fact that the group $WB_n$ is a homomorphic image of the group $VB_n$, it could be concluded from  Theorem 3 that the group $G(\beta_w)$ is an invariant of the link $L_w$.


\vspace{0.8cm}

\begin{center}
{\bf \S~6. Proof of  Theorem 3}
\end{center}

\vspace{0.5cm}

The proof of this theorem is divided into several steps. At first it will be shown that the group $G_{\widetilde{M}}(\beta)$ which is constructed using the representation $\widetilde{\varphi_M}$ is the invariant of closed braid $\widehat{\beta}$. It is true

\medskip

{\bf Theorem 4.}
{\it Let $\beta \in VB_n$ and $\beta' \in VB_m$ be two virtual braids such that theirs closures define the same link $L$, then $G_{\widetilde{M}}(\beta) \cong G_{\widetilde{M}}(\beta')$,
i.~e. the group $G_{\widetilde{M}}(\beta)$ is the invariant of the link
$L = \widehat{\beta}$ and it will be denoted by the  symbol $G_{\widetilde{M}}(L)$.}

\medskip

{\bf Proof.}
Kamada proved the analogue of Markov's theorem in the case of virtual braids \cite{SK}.
This theorem formulated in the algebraic terms by Kauffman and Lambropoulou \cite{K} is more convenient for the verification. According to this theorem two oriented virtual links are equivalent if and only if two corresponding virtual braids differ by braid relations in $VB_n$ and a finite sequence of the following moves or their inverses:

1) Virtual and real conjugation:
$$
\beta \to \rho_k \beta \rho_k,~~~ \beta \to \sigma_k \beta \sigma_k^{-1};
$$

2) Right virtual and real stabilization:
$$
\beta \to \beta \rho_n, ~~~ \beta \to \beta\sigma_n^{\pm1},
$$

3)  Algebraic right over/under threading:
$$
\beta \to \beta \sigma_n^{\pm1} \rho_{n-1} \sigma_n^{\mp1},
$$

4)  Algebraic left over/under threading:
$$
\beta \to \beta \rho_n \rho_{n-1} \sigma_{n-1}^{\mp1} \rho_n \sigma_{n-1}^{\pm1} \rho_{n-1} \rho_n,
$$
where $ \beta, \rho_k, \sigma_k \in VB_n, k = 1,\ldots, n-1$ and $ \rho_n, \sigma_n \in VB_{n+1}.$

It is needed to verify that groups obtained by moves 1) -- 4) are isomorphic to the group $G_{\widetilde{M}}(\beta)$.

To simplify notations we will write $\beta$ instead of $\widetilde{\varphi}_{M}(\beta)$ and write automorphisms to the right of the arguments.
We denote by symbol $G_1$ the group constructed using the braid $\beta$:
$$
G_1 = G_{\widetilde{M}}(\beta) = \langle y_1,y_2, \ldots,y_n, v_1, v_2, \ldots, v_n ~\| ~ y_i = y_i \beta,~~v_i = v_i \beta, ~i=1,2,\ldots,n \rangle.
$$
From now on we will not write commutativity relations of the elements $v_i$.

1) We consider the conjugacy of the braid $\beta$  in the group $VB_n$.
Conjugating by the generator  $\sigma_k^{-1}$ the following group is obtained
\[
G_2 = G_{\widetilde{M}}(\sigma_k\beta\sigma_k^{-1}) = \langle y_1, \ldots,y_n, v_1,  \ldots, v_n ~ \| ~ y_i = y_i (\sigma_k\beta\sigma_k^{-1}), ~v_i = v_i (\sigma_k\beta\sigma_k^{-1}), ~i=1,2,\ldots,n\rangle,
\]
where $k \in {1,2,\ldots,n-1}$. To prove that $G_2 \cong G_1$ defining relations of $G_2$ are rewritten in the form
\[
y_i\sigma_k = y_i(\sigma_k\beta), ~i = 1,2,\ldots,n.
\]
If $i \neq k, k+1$, due to the equality $y_i \sigma_k = y_i$, this relation is equivalent consequently
\[
y_i = y_i \beta
\]
which is the relation in $G_1$.  Thus only four relations are needed to consider:
\begin{equation}
\begin{split}
y_k \sigma_k = y_k(\sigma_k \beta), \hspace{1cm} y_{k+1}\sigma_k = y_{k+1}(\sigma_k \beta),\\
v_k \sigma_k = v_k(\sigma_k \beta), \hspace{1cm} v_{k+1}\sigma_k = v_{k+1}(\sigma_k \beta).  \label{a}
\end{split}
\end{equation}

By definition of $\widetilde{\varphi}_{M}$ relations (\ref{a}) are equivalent to
\begin{align*}
\begin{split}
y_k y_{k+1} y_k^{-1} &= (y_k y_{k+1} y_k^{-1}) \beta, \hspace{0,7cm} y_k = y_k \beta,\\
v_{k+1} &= v_{k+1} \beta, \hspace{2cm} v_k = v_k \beta.
\end{split}
\end{align*}
It is clear that the second, the third and  the fourth relations are relations in $G_1$. We rewrite the first relation:
$$
y_k y_{k+1} y_k^{-1} = (y_k \beta) (y_{k+1} \beta) (y_k^{-1} \beta)
$$
Using the second relation we obtain
$$
y_{k+1} = y_{k+1} \beta,
$$
i. e. it is a relation in $G_1$. Hence it has been proved that it is possible to go from the group $G_2$ to the group $G_1$ using the Tietze transformations.

We are considering the conjugation by the element $\rho_k$. In that case we have
\[
G_2 = G_{\widetilde{M}}(\rho_k \beta \rho_k) =
 \langle y_1, \ldots,y_n, v_1,  \ldots, v_n ~ \| ~ y_i = y_i (\rho_k \beta \rho_k),
~v_i = v_i (\rho_k \beta \rho_k), ~i=1,2,\ldots,n\rangle.
\]
We rewrite defining relations of $G_2$ in the form
\[
y_i\rho_k = y_i(\rho_k\beta), ~i = 1,2,\ldots,n.
\]
If $i \neq k, k+1$, then we obtain
\[
y_i = y_i \beta,
\]
due to $y_i \rho_k = y_i$. But it is the relation in $G_1$.
Hence only four relations have to be considered:
\begin{align*}
\begin{split}
y_k \rho_k = y_k(\rho_k \beta), \hspace{1cm} y_{k+1}\rho_k = y_{k+1}(\rho_k \beta),\\
v_k \rho_k = v_k(\rho_k \beta), \hspace{1cm} v_{k+1}\rho_k = v_{k+1}(\rho_k \beta).
\end{split}
\end{align*}

These relations are equivalent to the following ones
\begin{align*}
\begin{split}
y_{k+1}^{v_k^{-1}} &= y_{k+1}^{v_k^{-1}} \beta, \hspace{1cm} y_{k}^{v_{k+1}} = y_{k}^{v_{k+1}} \beta,\\
v_{k+1} &= v_{k+1} \beta, \hspace{1,6cm} v_k = v_k \beta.
\end{split}
\end{align*}
It is easy to see that the relation of $G_1$ is obtained.
Thus the fact has been proved that the set of relations in $G_2$ is equivalent to the set of relations in $G_1$.

2) We are considering the move of the braid $\beta  \in VB_n $ into the braid
$\beta \sigma_n^{-1} \in VB_{n+1} $. We are dealing with two groups $G_1$ and
$$
G_2 = G_{\widetilde{M}}(\beta\sigma_n^{-1}) = \langle y_1, \ldots,y_n, y_{n+1}, v_1,  \ldots, v_n, v_{n+1} ~ \| ~ y_i = y_i (\beta \sigma_n^{-1}), ~v_i = v_i (\beta \sigma_n^{-1}), ~i=1,2,\ldots,n+1\rangle
$$
and it is needed to prove that they are isomorphic.

We rewrite relations of $G_2$ in the form
\[
y_i\sigma_n = y_i\beta, ~i = 1,2,\ldots,n+1.
\]
Using the same reasoning as in 1) we have to consider only four relations:
\begin{equation}
\begin{split}
y_n \sigma_n = y_n \beta, \hspace{1cm} y_{n+1}\sigma_n = y_{n+1}\beta\\
v_n \sigma_n = v_n \beta, \hspace{1cm} v_{n+1}\sigma_n = v_{n+1}\beta. \label{r}
\end{split}
\end{equation}

By definition $\varphi_{\widetilde{M}}$ relations (\ref{r}) are equivalent to
\begin{align*}
\begin{split}
y_n y_{n+1} y_n^{-1} &= y_n \beta, \hspace{1cm} y_n = y_{n+1},\\
v_{n+1} &= v_n \beta, \hspace{1cm} v_n = v_{n+1}.
\end{split}
\end{align*}
Using the second relation the first one is rewritten in the form
$$
y_n = y_n \beta
$$
which is the relation in $G_1$. Using relations $y_{n+1} = y_n$ and $v_{n+1} = v_n$ we can exclude $y_{n+1}$ and $v_{n+1}$ from the generator set of the group $G_2$.
Thus it has been proved that the group $G_2$ is isomorphic to the group $G_1$.

The move of the braid $\beta  \in VB_{n} $ into the braid $\beta \sigma_n \in VB_{n+1} $ are considered analogously.

Now we consider the move of the braid $\beta  \in VB_n $ into the braid $\beta \rho_n \in VB_{n+1}$.
Two groups are obtained, $G_1$ and
$$
G_2 = G_{\widehat{M}}(\beta\rho_n) = \langle y_1, \ldots,y_n, y_{n+1}, v_1,  \ldots, v_n, v_{n+1} ~ \| ~ y_i = y_i (\beta \rho_n), ~v_i = v_i (\beta \rho_n), ~i=1,2,\ldots,n+1\rangle.
$$

For  $ i = n$ and $i = n+1$ we have the following relations in $G_2$:
\begin{align*}
\begin{split}
y_n \rho_n = y_n \beta, \hspace{1cm} y_{n+1}\rho_n = y_{n+1}\beta\\
v_n \rho_n = v_n \beta, \hspace{1cm} v_{n+1}\rho_n = v_{n+1}\beta
\end{split}
\end{align*}
which are equivalent to relations
\begin{align*}
\begin{split}
y_{n+1}^{v_n^{-1}} &= y_{n} \beta, \hspace{1cm} y_{n}^{v_{n+1}} = y_{n+1} \beta = y_{n+1},\\
v_{n+1} &= v_n \beta, \hspace{1cm} v_n = v_{n+1}.
\end{split}
\end{align*}

Substituting $y_{n+1}$ from the second equation into the first relation, we obtain
$$
y_n = y_n \beta
$$
and that is the relation in $G_1$. Using the second and the fourth relations variables $y_{n+1}$ and $v_{n+1}$ could be excluded from the set of generators  of the group $G_2$.
Hence we have proved that the group $G_2$ is isomorphic to the group $G_1$.

3) We are considering the algebraic right over threading,
i. e. the move of the braid $\beta \in VB_n $ into the braid
$\beta \sigma_n \rho_{n-1} \sigma_n^{-1} \in VB_{n+1}$. In this case we have $G_1$ and
\begin{align*}
\begin{split}
G_2 = G_{\widetilde{M}}(\beta \sigma_n \rho_{n-1} \sigma_n^{-1}) =& \langle y_1,y_2,\ldots,y_n, y_{n+1}, v_1,  \ldots, v_n, v_{n+1} ~ \| ~ \\
 y_i = y_i (\beta \sigma_n &\rho_{n-1} \sigma_n^{-1}),~~v_i = v_i (\beta \sigma_n \rho_{n-1} \sigma_n^{-1}), ~~i=1,2,\ldots,{n+1}\rangle.
\end{split}
\end{align*}
To prove that $G_1 \cong G_2$ we rewrite relations of $G_2$ in the form
\[
y_i \sigma_n \rho_{n-1} \sigma_n^{-1} =
 y_i\beta,~~ v_i \sigma_n \rho_{n-1} \sigma_n^{-1} = v_i\beta, ~i = 1,2,\ldots,n+1.
\]
If $i = 1, 2,\ldots, n-2$ then we have
\[
y_i = y_i \beta, ~~v_i = v_i \beta.
\]
These are relations in $G_1$.  Hence we have to consider only six relations:
\begin{equation}
\begin{split}
y_{n-1}(\sigma_n \rho_{n-1} \sigma_n^{-1}) =
 y_{n-1}\beta, ~~y_n (\sigma_n \rho_{n-1} \sigma_n^{-1}) =
   y_n \beta, ~~y_{n+1}(\sigma_n \rho_{n-1} \sigma_n^{-1}) = y_{n+1}\beta, \\
v_{n-1}(\sigma_n \rho_{n-1} \sigma_n^{-1}) =
 v_{n-1}\beta, ~~v_n (\sigma_n \rho_{n-1} \sigma_n^{-1}) =
   v_n \beta, ~~v_{n+1}(\sigma_n \rho_{n-1} \sigma_n^{-1}) = v_{n+1}\beta. \label{u}
\end{split}
\end{equation}

By definition $\widetilde{\varphi}_{M}$ the braid $\sigma_n \rho_{n-1} \sigma_n^{-1} $
acts on generators $y_{n-1}, y_n, y_{n+1},v_{n-1}, v_n, v_{n+1}$ in the following way:
\begin{align*}
\begin{split}
\sigma_n \rho_{n-1} \sigma_n^{-1} :\left\{%
\begin{array}{ll}
    y_{n-1} \xrightarrow {\sigma_n} y_{n-1} \xrightarrow {\rho_{n-1}} y_n^{v_{n-1}^{-1}} \xrightarrow {\sigma_n^{-1}} y_{n+1}^{v_{n-1}^{-1}}, \\
    y_{n} \xrightarrow {\sigma_n}  y_n y_{n+1} y_n^{-1} \xrightarrow {\rho_{n-1}} y_{n-1}^{v_n} y_{n+1} y_{n-1}^{-v_n}  \xrightarrow {\sigma_n^{-1}} y_{n-1}^{v_{n+1}} y_{n+1}^{-1} y_n y_{n+1} y_{n-1}^{-v_{n+1}} , \\
    y_{n+1} \xrightarrow {\sigma_n} y_n \xrightarrow {\rho_{n-1}} y_{n-1}^{v_n} \xrightarrow {\sigma_n^{-1}}  y_{n-1}^{v_{n+1}}.\\
\end{array}%
\right.
\end{split}
\end{align*}
Analogously
\begin{align*}
\sigma_n \rho_{n-1} \sigma_n^{-1} :\left\{%
\begin{array}{ll}
    v_{n-1} \xrightarrow {\sigma_n \rho_{n-1} \sigma_n^{-1}} v_{n+1}, \\
    v_{n} \xrightarrow {\sigma_n \rho_{n-1} \sigma_n^{-1}}  v_n , \\
    v_{n+1} \xrightarrow {\sigma_n \rho_{n-1} \sigma_n^{-1}} v_{n-1}.\\
\end{array}%
\right.
\end{align*}

Therefore, relations (\ref{u}) for the representation $\widetilde{\varphi}_{M}$ are equivalent to
$$
y_{n+1}^{v_{n-1}^{-1}} = y_{n-1}\beta,~~ y_{n-1}^{v_{n+1}} y_{n+1}^{-1} y_n y_{n+1} y_{n-1}^{-v_{n+1}} = y_{n}\beta, ~~ y_{n-1}^{v_{n+1}} = y_{n+1}\beta,
$$
$$
v_{n+1} = v_{n-1}\beta, ~~ v_n = v_{n}\beta,~~ v_{n-1} = v_{n+1}\beta.
$$
The braid $\beta$ acts on generators $y_{n+1}$ and $v_{n+1}$
trivially since  $\beta \in VB_n$, i. e. $y_{n+1}\beta = y_{n+1},~ v_{n+1}\beta = v_{n+1}$.
Consequently, the third and the sixth relations are equivalent to
$$
y_{n-1}^{v_{n+1}} = y_{n+1}, ~~ v_{n-1} = v_{n+1},
$$
correspondingly.
Substituting the third relation into the first one and taking into account the sixth relation, we obtain
$$
y_{n-1} = y_{n-1}\beta,
$$
which is the relation in $G_1$.
Using the third relation we obtain that the second relation is equivalent to
$$
y_{n} = y_{n}\beta,
$$
which is the relation in $G_1$. Besides $y_{n+1}$ and $v_{n+1}$ could be excluded from the generator set of $G_2$. Therefore, it has been proved that the set of relations of $G_2$ is equivalent to the set of relations of $G_1$.

4) We are considering the algebraic left over threading, i.~e. the move of the braid $\beta \in VB_n $ into the braid
 $\beta \rho_n \rho_{n-1} \sigma_{n-1}^{-1} \rho_{n} \sigma_{n-1} \rho_{n-1} \rho_n \in VB_{n+1} $. In this case
\begin{align*}
\begin{split}
G_2 = G(\beta \sigma_n \rho_{n-1} \sigma_n^{-1})& = \langle y_1,y_2,\ldots,y_n, y_{n+1}, v_1,\ldots,v_n, v_{n+1} ~ \| \\
~ y_i= y_i (\beta \rho_n \rho_{n-1} \sigma_{n-1}^{-1} \rho_{n} \sigma_{n-1} \rho_{n-1} \rho_n)&,~
v_i = v_i (\beta \rho_n \rho_{n-1} \sigma_{n-1}^{-1} \rho_{n} \sigma_{n-1} \rho_{n-1} \rho_n), ~i=1,2,\ldots,{n+1}\rangle.
\end{split}
\end{align*}
To prove that $G_1 \cong G_2$ we rewrite relations of  $G_2$ in the form
\[
y_i \rho_n \rho_{n-1} \sigma_{n-1}^{-1} \rho_{n} \sigma_{n-1} \rho_{n-1} \rho_n = y_i\beta, ~i = 1,2,\ldots,n+1.
\]
If $i = 1, 2,\ldots, n-2$ we have
\[
y_i = y_i \beta
\]
which is the relation in $G_1$.  Hence it has to be considered only six relations:
\begin{equation}
\begin{split}
y_{n-1}(\rho_n \rho_{n-1} \sigma_{n-1}^{-1} \rho_{n} \sigma_{n-1} \rho_{n-1} \rho_n) &= y_{n-1}\beta, \\
y_n (\rho_n \rho_{n-1} \sigma_{n-1}^{-1} \rho_{n} \sigma_{n-1} \rho_{n-1} \rho_n) &= y_n \beta, \\
y_{n+1}(\rho_n \rho_{n-1} \sigma_{n-1}^{-1} \rho_{n} \sigma_{n-1} \rho_{n-1} \rho_n) &= y_{n+1}\beta,\\
v_{n-1}(\rho_n \rho_{n-1} \sigma_{n-1}^{-1} \rho_{n} \sigma_{n-1} \rho_{n-1} \rho_n) &= v_{n-1}\beta, \\
v_n (\rho_n \rho_{n-1} \sigma_{n-1}^{-1} \rho_{n} \sigma_{n-1} \rho_{n-1} \rho_n) &= v_n \beta, \\
v_{n+1}(\rho_n \rho_{n-1} \sigma_{n-1}^{-1} \rho_{n} \sigma_{n-1} \rho_{n-1} \rho_n) &= v_{n+1}\beta. \label{l}
\end{split}
\end{equation}

Denote the braid $\rho_n \rho_{n-1} \sigma_{n-1}^{-1} \rho_{n} \sigma_{n-1} \rho_{n-1} \rho_n$ by $b$. It is easy to see that
$b$
acts on generators $y_{n-1}, y_n, y_{n+1}$, $v_{n-1}, v_n, v_{n+1}$ in the following way:
\begin{align*}
b :\left\{%
\begin{array}{ll}
    y_{n-1} \longrightarrow (y_n^{-1} y_{n+1}^{v_n^{-1} v_{n-1}^{-1}} y_{n-1}^{v_{n+1}} y_{n+1}^{-v_n^{-1} v_{n-1}^{-1}} y_n)^{v_{n}^{-1}},\\
    y_{n} \longrightarrow  y_{n+1}^{v_n^{-1} v_{n-1}^{-1}}, \\
    y_{n+1} \longrightarrow y_{n}^{v_{n+1} v_{n-1}}.\\
\end{array}%
\right.
\end{align*}
\begin{align*}
b :\left\{%
\begin{array}{ll}
    v_{n-1} \longrightarrow v_{n-1},\\
    v_{n} \longrightarrow  v_{n+1}, \\
    v_{n+1} \longrightarrow v_{n}.\\
\end{array}%
\right.
\end{align*}

Thus relations (\ref{l}) are equivalent to
$$
(y_n^{-1} y_{n+1}^{v_n^{-1} v_{n-1}^{-1}} y_{n-1}^{v_{n+1}} y_{n+1}^{-v_n^{-1} v_{n-1}^{-1}} y_n)^{v_{n}^{-1}} = y_{n-1}\beta, ~~
y_{n+1}^{v_n^{-1} v_{n-1}^{-1}} = y_{n}\beta, ~~ y_{n}^{v_{n+1} v_{n-1}} = y_{n+1}\beta,
$$
$$
v_{n-1} = v_{n-1}\beta, ~~ v_{n+1} = v_{n}\beta, ~~ v_{n} = v_{n+1}\beta.
$$
Since $\beta \in VB_n$ it acts on generators $y_{n+1}$ and $v_{n+1}$ trivially. Therefore, the third relation is equivalent to
$$
y_{n}^{v_{n+1} v_{n-1}} = y_{n+1}
$$
and the sixth relation is equivalent to
$$
v_{n} = v_{n+1}.
$$
Obtained equations are substituted in the first and the second relations and it results in
$$
y_{n-1} = y_{n-1}\beta, ~~y_{n} = y_{n}\beta,
$$
which are the relations in $G_1$.
The generator  $y_{n+1}$ and $v_{n+1}$ could be excluded from the generating set of $G_2$.
Therefore, the fact has been proved that the set of relations of $G_2$ is equivalent to the set of relations of $G_1$.

The theorem is proved.
\medskip

We proved that the representation $\varphi_M$ is equivalent to the representation $\widetilde{\varphi}_M$. Now we show that the group $G_M(\beta)$ is isomorphic to the group $G_{\widehat{M}}(\beta)$. According to the proved theorem, the group $G_M(\beta)$ is an invariant of the link $\widehat{\beta}$ and this fact leads to the proof of Theorem 3.

Let a group $G$ have the presentation
$G=\langle  X \, \| \,\mathcal{ R}  \rangle$,
where $X=\left\{ x_1,\ldots,x_n \right\}$ is the generating set and
$\mathcal{R}$ is the set of defining relations.
If in the group $G$ another generating set $Y=\left\{ y_1,\ldots,y_m \right\}$ is chosen,
then the group $G$ has the presentation  $G=\langle  Y \, \| \, \widetilde{\mathcal{R}}  \rangle$
which connected with the initial presentation by the sequence of Tietze transformations.
From now on we assume the generating sets $X$, $Y$ and their corresponding sets of defining relations are fixed.

If  $\varphi$ is some automorphism of the group $G$, then a group could be defined using the generating set $X$
$$
G(X,\varphi)=\langle  X \, \| \, \mathcal{R},\,\, \varphi(x_i) = x_i,\, i=1,\ldots, n  \rangle.
$$
Analogously using the generating set $Y$ and the automorphism  $\varphi$, a group could be defined
$$
G(Y,\varphi)=\langle  Y \, \| \, \widetilde{\mathcal{R}},\,\, \varphi(y_j) = y_j,\, j=1,\ldots, m  \rangle.
$$
It is true

{\bf Theorem 5.}
{\it
Groups $G(X,\varphi)$  and $G(Y,\varphi)$ are isomorphic.
}

{\bf Proof.}
Let the automorphism $\varphi$ act on generating sets $X$ and $Y$
in the following way
$$
\varphi(x_i) = w_i(X), \quad \varphi(y_j) = u_j(Y),
$$
where $w_i(X)$ and $u_j(Y)$ are words in the terms of generators $X$ and $Y$ as well as their inverses consequently for
$i=1,\ldots, n$, $j=1,\ldots, m$.
We fix expressions
$$
x_i=b_i(Y), \quad y_j=a_j(X)
$$
to express generators of $X$ through $Y$ and, visa versa, generators $Y$ through $X$.
To simplify the notation indices  $i$ and $j$ are dropped out so there are obtained systems of equations
$$
\varphi(X) = W(X), \quad \varphi(Y) = U(Y), \quad
X=B(Y), \quad Y=A(X).
$$
Since the automorphism $\varphi$ of the group $G$ that acts on generators of $X$ and $Y$ is the same, the following systems of equations
$$
A(W(B(Y)))=U(Y) \Leftrightarrow
                 A(W(X))=U(A(X)) \Leftrightarrow
                                   W(B(Y))=B(U(Y))
$$
are equivalent in the group $G$.

The proof of the theorem statement is to apply sequentially the Tietze transformations to the presentation  of the group $G(X,\varphi)$.
Thus we consider the group $G(X,\varphi)$.
By definition we have
$$
G(X,\varphi)=\langle  X \, \| \, \mathcal{R},\,\, X=W(X)  \rangle.
$$
We add elements $Y$ into the generating set and the equations $Y=A(X)$ into the set of relations
$$
G(X,\varphi)=\langle  X, \, Y\, \| \, \mathcal{R},\,\, X=W(X),\,\, Y=A(X)  \rangle.
$$
Relations $Y=A(X)$ are equivalent to  $X=B(Y)$ with the module of relations $\mathcal{R}$.
Thus
$$
\langle  X, \, Y\, \| \, \mathcal{R},\,\, X=W(X),\,\, Y=A(X)  \rangle \simeq
\langle  X, \, Y\, \| \, \mathcal{R},\,\, X=W(X),\,\, X=B(Y)  \rangle.
$$
Excluding generators $X$ and relations $X=B(Y)$ we obtain
$$
\langle  X, \, Y\, \| \, \mathcal{R},\,\, X=W(X),\,\, X=B(Y)  \rangle \simeq
\langle  Y\, \| \, \widetilde{\mathcal{R}},\,\, B(Y)=W(B(Y))  \rangle.
$$
Since $W(B(Y))=B(U(Y))$ we have
$$
\langle  Y\, \| \, \widetilde{\mathcal{R}},\,\, B(Y)=W(B(Y))  \rangle \simeq
\langle  Y\, \| \, \widetilde{\mathcal{R}},\,\, B(Y)=B(U(Y))  \rangle.
$$
Taking into account that in the group $G$ the equations $A(B(Y))=Y$
and $B(A(X))=X$ hold, it could be seen that systems of relations $B(Y)=B(Z)$ and $Y=Z$ are equivalent in the group $G$ (i.~e. they are equivalent with the module of relations $\widetilde{\mathcal{R}}$). Consequently
$$
 \langle  Y\, \| \, \widetilde{\mathcal{R}},\,\, B(Y)=B(U(Y))  \rangle \simeq
 \langle  Y\, \| \, \widetilde{\mathcal{R}},\,\, Y=U(Y)  \rangle=
 G(Y,\varphi).
$$
Therefore,
$$
 G(X,\varphi) \simeq G(Y,\varphi).
$$
The theorem is proved.

\vspace{0.5cm}

As a corollary the fact is obtained that groups $G_M(L)$ and $G_{\widetilde{M}}(L)$ are isomorphic. Consequently Theorem 3 is proved.

\vspace{0.8cm}

\begin{center}
{\bf \S~7. Properties of link groups}
\end{center}

\vspace{0.5cm}

Directly from the definition of the group $G_M$ it could be derived

{\bf Proposition 6.} {\it If $L$ is $d$-component link which is
presented by a closure of  a $n$-strand braid, then in the group $G_M(L)$ there are $d$ different generators among $u_1, u_2, \ldots, u_n$.
In particular, provided $L$ a knot, $u_1 = u_2 = \ldots = u_n$.}

\medskip

As shown earlier the representation $\widetilde{\varphi}_M$ is an extension of the Artin representation. As a corollary one can prove

\medskip

{\bf Proposition 7.} {\it Given  a classical $d$-component link $L$,
there is the following isomorphism of groups}
$$
G_{\widetilde{M}} (L) = G(L) * \mathbb{Z}^d,~\mbox{where}~G(L) = \pi_1(\mathbb{S}^3 \setminus L).
$$

\begin{figure}[h]
\noindent\centering{\includegraphics[height=0.4 \textwidth]{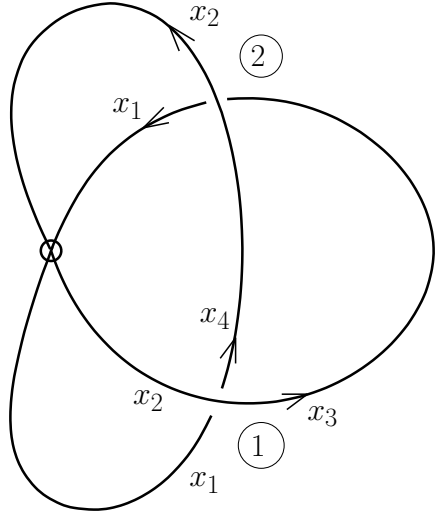}}
\caption{Virtual trefoil}
\label{vtrefoild}
\end{figure}
{\bf Example.} Given $T_v$ the virtual trefoil (see Fig. \ref{vtrefoild}),
the generalized Alexander group is generated by elements $x_1$, $x_2$, $x_3$, $x_4$, $u$, $v$ and it has in the first crossing the relations
$$
x_1 x_2^u = x_3 x_4^u,~~x_2 = x_3^v,
$$
and in the second crossing the relations
$$
x_3 x_4^u = x_2 x_1^u,~~x_4 = x_2^v.
$$
Excluding generators $x_3$ and $x_4$, we obtain the presentation
$$
\widetilde{A}_{T_v} =
\langle x_1, x_2, u, v ~||~[u,v]=1,~~x_1 x_2^{u} = x_2^{v^{-1}} x_2^{v u},~~x_2^{v^{-1}} x_2^{v u} = x_2 x_1^{u} \rangle.
$$
Excluding the generator $x_1$, we have
$$
\widetilde{A}_{T_v} =
\langle x_2, u, v ~||~[u,v]=1,~~x_2 x_2^{v^{-1} u} x_2^{vu^2} = x_2^{v^{-1}} x_2^{v u} x_2^{u^2} \rangle.
$$
If a generator $x_2$ replaced by a new one $b$ which equals to $x_2 u^{-1}$, then the presentation is of the form
$$
\widetilde{A}_{T_v} =
\langle b, u, v ~||~[u,v]=1,~~b (b^{v^{-1}} b^{v}) = (b^{v^{-1}} b^{v}) b \rangle = \langle b, u, v ~||~[u,v]=1,~~[b^{v^{-1}} b^{v}, b] = 1 \rangle.
$$


On the other hand, the virtual trefoil could be presented as a classical knot in the thickened torus. The fundamental group of the knot complement in the thickened torus was found in \cite{CSW}
and it has the following presentation
$$
G = \langle a, x, y ~||~[x,y]=1,~~(a^{x} a^{y}) a^{x y} = a (a^{x} a^{y}) \rangle.
$$
A natural question arises about the isomorphism of groups $\widetilde{A}_{T_v}$ and $G$. The answer on this question is given by
\medskip

{\bf Proposition 8. }
{\it
Groups $\widetilde{A}_{T_v}$ and $G$ are not isomorphic.
}

\medskip

{\bf Proof.} For brevity we denote $H=\widetilde{A}_{T_v}$.
We show that the factor-groups $H/\gamma_3 H$ and $G/\gamma_3 G$ of these groups by the third term of lower central series are not isomorphic to each other.
We obtain
$$
b b^{v^{-1}} b^{v} =b^2[b,v^{-1}]b[b,v]\equiv b^3 \, (\mathrm{mod}\, \gamma_3 H),
$$
$$
b^{v^{-1}} b^{v} b =b[b,v^{-1}]b[b,v]b \equiv b^3 \, (\mathrm{mod}\, \gamma_3 H).
$$
Hence
$$
H/\gamma_3 H =
\langle b, u, v ~||~[u,v]=1,~~\gamma_3 H=1 \rangle.
$$
In particular, $\gamma_2 H/\gamma_3 H \cong \mathbb{Z} \times \mathbb{Z}$
(it is generated by images of the commutators $[b,u]$, $[b,v]$).

On the other hand
$$
a^{x} a^{y} a^{x y} =a[a,x]a[a,y]a[a,xy]\equiv a^3 [a,xy]^2 \, (\mathrm{mod}\, \gamma_3 G),
$$
$$
a a^{x} a^{y} =a^2[a,x]a[a,y]\equiv a^3 [a,xy] \, (\mathrm{mod}\, \gamma_3 G).
$$
Therefore,
$$
G/\gamma_3 G =
\langle a, x, y ~||~[x,y]=1,~~[a,xy]=1,~~\gamma_3 G=1 \rangle.
$$
In particular, $\gamma_2 G/\gamma_3 G \cong \mathbb{Z}$
(it is generated by images of the commutator $[a,x]$).

As factor-groups
$\gamma_2 H/\gamma_3 H$ and
$\gamma_2 G/\gamma_3 G$ are not isomorphic so the correspondent groups
$H$ and $G$ are not isomorphic.
The proposition is proved.

We show that the group $\widetilde{A}_{T_v}$ is a
HNN-extension of the one-relator group.

{\bf Proposition 9. }
{\it
The group
$$
\widetilde{A}_{T_v} =
\langle b, u, v ~||~[u,v]=1,~~b b^{v^{-1}} b^{v} = b^{v^{-1}} b^{v} b \rangle
$$
is a $HNN$-extension.
}

{\bf Proof.}
Using Tietze transformations for the presentation of $\widetilde{A}_{T_v}$,  we obtain
$$
\widetilde{A}_{T_v} =
\langle b, u, v ~||~[u,v]=1,~~b b^{v^{-1}} b^{v} = b^{v^{-1}} b^{v} b \rangle=
$$
$$
 =\langle b, u, v ~||~u^v=u,~~ [b^{v^{-1}} b^{v},b] = 1 \rangle=
$$
$$
 =\langle b, u, v, z, t ~||~u^v=u,~~z=b^{v^{-1}},~~t=b^v,~~ [zt,b] = 1 \rangle=
$$
$$
 =\langle b, u, v, z, t ~||~u^v=u,~~z^v=b,~~b^v=t,~~ [zt,b] = 1 \rangle.
$$

We consider the group
$$
 G=\langle b, u, z, t ~||~ [zt,b] = 1 \rangle
$$
and its subgroups
$$
 A=\langle b, u, z  \rangle, \quad B=\langle b, u, t  \rangle.
$$
It is sufficient to show that the mapping
$$
 u \mapsto u, \quad z \mapsto b, \quad b \mapsto t
$$
is an isomorphism of the subgroup $A$ onto the subgroup $B$. For that purpose it is sufficient to show that
$A$ and $B$ are free subgroups of the rank 3
with free generators $u, b, z$ and $u, b, t $ respectively.
Taking into account that
$$
 G=\langle b, u, z, t ~||~ [zt,b] = 1 \rangle=\langle u \rangle \ast \widetilde{G},
$$
where
$$
 \widetilde{G}=\langle b, z, t ~||~ [zt,b] = 1 \rangle,
$$
it is sufficient to show that subgroups
$$
 \widetilde{A}=\langle b, z  \rangle, \quad \widetilde{B}=\langle b, t  \rangle
$$
are free groups of the rank 2 in the group $\widetilde{G}$.
We consider the homomorphic image of the group $\widetilde{G}$ which is obtained by adding relations $t=z^{-1}$:
$$
 \langle b, z, t ~||~ [zt,b] = 1,~~t=z^{-1} \rangle.
$$
Excluding the generator $t$, the free group $\langle b, z  \rangle$ is obtained which is covered by the image of the subgroup $\widetilde{A}$.
And excluding the generator $z$, the free group $\langle t, z  \rangle$ is obtained which is covered by the image of the subgroup $\widetilde{B}$.
Therefore, $\widetilde{A}$ and  $\widetilde{B}$ are free subgroups of the rank 2 with free generators $b, z$ and $b, t $ respectively.
It leads to the fact that $A$ and  $B$ are  free subgroups of the rank 3
with free generators $u, b, z$ and $u, b, t $ respectively.

Hence
$$
\widetilde{A}_{T_v}
 =\langle b, u, v, z, t ~||~u^v=u,~~z^v=b,~~b^v=t,~~ [zt,b] = 1 \rangle
$$
is a HNN-extension of the group
$$
 G=\langle b, u, z, t ~||~ [zt,b] = 1 \rangle
$$
with the stable letter $v$ and associated subgroups
$$
 A=\langle b, u, z  \rangle, \quad B=\langle b, u, t  \rangle.
$$

The proposition is proved.

\end{document}